\documentclass[10pt]{amsart}
\input{epsf}

\usepackage{amssymb,latexsym}
\topskip  \title[On multi-avoidance of generalized patterns]{On multi-avoidance of generalized patterns}
\author{Sergey Kitaev and Toufik Mansour}
\address{Matematik, Chalmers tekniska h\"ogskola och G\"oteborgs
universitet, 412~96  G\"oteborg, Sweden}
\email{kitaev@math.chalmers.se} 
\address{Department of mathematics, Chalmers university of
technology, 412~96 G\"oteborg, Sweden}
\email{toufik@math.chalmers.se} 

\newtheorem{prop}{Proposition}
\newtheorem{lemma}[prop]{Lemma}

\newtheorem{thm}[prop]{Theorem}
\theoremstyle{definition}

\def\mn{{\mbox{-}}}

\begin{document}
\begin{abstract}
In~\cite{Kit1} Kitaev discussed simultaneous avoidance of two
3-patterns with no internal dashes, that is, where the patterns
correspond to contiguous subwords in a permutation. In three
essentially different cases, the numbers of such $n$-permutations
are $2^{n-1}$, the number of involutions in $\mathcal{S}_n$, and
$2E_n$, where $E_n$ is the $n$-th Euler number. In this paper we
give recurrence relations for the remaining three essentially
different cases.

To complete the descriptions in \cite{Kit3} and \cite{KitMans}, we
consider avoidance of a pattern of the form $x\mn y\mn z$ (a
classical 3-pattern) and beginning or ending with an increasing or
decreasing pattern. Moreover, we generalize this problem: we
demand that a permutation must avoid a 3-pattern, begin with a
certain pattern and end with a certain pattern simultaneously. We
find the number of such permutations in case of avoiding an
arbitrary generalized 3-pattern and beginning and ending with
increasing or decreasing patterns.

\end{abstract}

\maketitle

\thispagestyle{empty}

\section{Introduction and Background}
{\bf Permutation patterns:} All permutations in this paper are
written as words $\pi=a_1 a_2\ldots a_n$, where the $a_i$ consist
of all the integers $1,2,\ldots,n$. Let $\alpha\in S_n$ and
$\tau\in S_k$ be two permutations. We say that $\alpha$ {\em
contains} $\tau$ if there exists a subsequence $1\leq
i_1<i_2<\cdots<i_k\leq n$ such that
$(\alpha_{i_1},\dots,\alpha_{i_k})$ is order-isomorphic to $\tau$;
in such a context $\tau$ is usually called a {\em pattern}. We say
that $\alpha$ {\em avoids} $\tau$, or is $\tau$-{\em avoiding}, if
such a subsequence does not exist. The set of all $\tau$-avoiding
permutations in $S_n$ is denoted by $S_n(\tau)$. For an arbitrary
finite collection of patterns $T$, we say that $\alpha$ avoids $T$
if $\alpha$ avoids any $\tau\in T$; the corresponding subset of
$S_n$ is denoted by $S_n(T)$.

While the case of permutations avoiding a single pattern has
attracted much attention, the case of multiple pattern avoidance
remains less investigated. In particular, it is natural, as the
next step, to consider permutations avoiding pairs of patterns
$\tau_1$, $\tau_2$. This problem was solved completely for
$\tau_1,\tau_2\in S_3$ (see \cite{SchSim}), for $\tau_1\in S_3$
and $\tau_2\in S_4$ (see \cite{W}), and for $\tau_1,\tau_2\in S_4$
(see \cite{B1,Km} and references therein). Several recent papers
\cite{CW,MV1,Kr,MV3,MV2} deal with the case $\tau_1\in S_3$,
$\tau_2\in S_k$ for various pairs $\tau_1,\tau_2$.

{\bf Generalized permutation patterns:} In \cite{BabStein} Babson
and Steingr\'{\i}msson introduced {\em generalized permutation
patterns {\rm(}GPs{\rm)}} where two adjacent letters in a pattern
may be required to be adjacent in the permutation. Such an
adjacency requirement is indicated by the absence of a dash
between the corresponding letters in the pattern. For example, the
permutation $\pi=516423$ has only one occurrence of the pattern
$2\mn31$, namely the subword 564, but the pattern $2\mn3\mn1$
occurs also in the subwords 562 and 563. Note that a classical
pattern should, in our notation, have dashes at the beginning and
end. Since most of the patterns considered in this paper satisfy
this, we suppress these dashes from the notation. Thus, a pattern
with no dashes corresponds to a contiguous subword anywhere in a
permutation. The motivation for introducing these patterns was the
study of Mahonian statistics. A number of results on GPs were
obtained by Claesson, Kitaev and Mansour. See for example
\cite{Claes}, \cite{Kit1, Kit2, Kit3} and \cite{Mans1, Mans2,
Mans3}.

As in \cite{SchSim}, dealing with the classical patterns, one can
consider the case when permutations have to avoid two or more
generalized patterns simultaneously. A complete solution for the
number of permutations avoiding a pair of 3-patterns of type (1,2)
or (2,1), that is the patterns having one internal dash, is given
in~\cite{ClaesMans1}. In~\cite{Kit1} Kitaev discussed simultaneous
avoidance of two 3-patterns with no internal dashes, that is,
where the patterns correspond to contiguous subwords in a
permutation. In three essentially different cases, the numbers of
such $n$-permutations are $2^{n-1}$, the number of involutions in
$\mathcal{S}_n$, and $2E_n$, where $E_n$ is the $n$-th Euler
number. The remaining cases are avoidance of $123$ and $231$,
$213$ and $231$, $132$ and $213$. In
Section~\ref{avoidanceOftwo3patterns} we give recurrence relations
for these cases.

In Section~4, we consider avoidance of a pattern $x\mn y\mn z$,
and beginning or ending with increasing or decreasing pattern.
This completes the results made in~\cite{KitMans}, which concerns
the number of permutations that avoid a generalized 3-pattern and
begin or end with an increasing or decreasing pattern.

In Sections~5--8, we give enumeration for the number of
permutations that avoid a generalazed 3-pattern, begin {\em and}
end with increasing or decreasing patterns. We record our results
in terms of either {\em generating functions}, or {\em exponential
generating functions}, or formulas for the numbers appeared.

In Section~\ref{section9}, we discuss possible directions of
generalization of the results from Sections~5--8.
\section{Preliminaries}\label{Preliminaries}

The {\em reverse} $R(\pi)$ of a permutation $\pi=a_1a_2 \ldots
a_n$ is the permutation $a_n \ldots a_2a_1$. The {\em complement}
$C(\pi)$ is the permutation $b_1b_2 \ldots b_n$ where
$b_i=n+1-a_i$. Also, $R \circ C$ is the composition of $R$ and
$C$. For example, $R(13254)=45231$, $C(13254)=53412$ and $R \circ
C(13254)=21435$. We call these bijections of $S_n$ to itself {\em
trivial}, and it is easy to see that for any pattern $p$ the
number $A_p(n)$ of permutations avoiding the pattern $p$ is the
same as for the patterns $R(p)$, $C(p)$ and $R \circ C(p)$. For
example, the number of permutations that avoid the pattern 132 is
the same as the number of permutations that avoid the pattern 231.
This property holds for sets of patterns as well. If we apply one
of the trivial bijections to all patterns of a set $G$, then we
get a set $G^{\prime}$ for which $A_{G^{\prime}}(n)$ is equal to
$A_G(n)$. For example, the number of permutations avoiding $\{
123, 132 \}$ equals the number of those avoiding  $\{ 321, 312\}$
because the second set is obtained from the first one by
complementing each pattern.

In this paper we denote the $n$th Catalan number by $C_n$; the
generating function for these numbers by $C(x)$; the $n$th Bell
number by $B_n$.

Also, $N_{p}^{q}(n)$ denotes the number of permutations that avoid
the pattern $p$ and begin with the pattern $q$; $G_{p}^{q}(x)$
(respectively, $E_{p}^{q}(x)$) denotes the ordinary (respectively,
exponential) generating function for the number of such
permutations. Besides, $N_{p}^{q, r}(n)$ denotes the number of
permutations that avoid the pattern $p$, begin with the pattern
$q$ and end with the pattern $r$; $G_{p}^{q, r}(x)$ (respectively,
$E_{p}^{q, r}(x)$) denotes the ordinary (respectively,
exponential) generating function for the number of such
permutations.

Recall the following properties of $C(x)$:
\begin{equation}
C(x) = \frac{1-\sqrt{1-4x}}{2x} = \frac{1}{1-xC(x)}.
\label{catalan}
\end{equation}
\section{Simultaneous avoidance of two 3-patterns with no dashes}\label{avoidanceOftwo3patterns}

\subsection{Avoidance of patterns 123 and 231 simultaneously}
We first consider the avoidance of the patterns 123 and
231 simultaneously.

Let $a(n;i_1,i_2,\dots,i_m)$ denote the number of permutations
$\pi\in S_n(123,231)$ such that $\pi_1\pi_2\dots\pi_m=i_1i_2\dots
i_m$ and let $a(n)=|S_n(123,231)|$.

By the definitions, we get that $a(n)=\sum_{j=1}^n a(n;j)$ and
$a(n;n)=a(n-1)$. Hence
\begin{equation}
a(n)=a(n-1)+a(n;1)+a(n;2)+\cdots+a(n;n-1).
\label{eq1}
\end{equation}

Also, by the definitions, for all $1\leq i\leq n-1$, we get
\begin{equation}
a(n;i)=\sum_{j=1}^{i-1}a(n;i,j)+\sum_{j=i+1}^na(n;i,j).
\label{eq2}
\end{equation}

Suppose $\pi\in S_n(123,231)$ is such that $\pi_1=i$ and
$\pi_2=j$. If $i>j$ then there is no occurrence of the pattern
$123$ or $231$ that contains $\pi_1$, so $a(n;i,j)=a(n-1;j)$. If
$i<j$ then since $\pi$ avoids $123$ and $231$, we get that
$i<\pi_3<j$, and thus in this case
$a(n;i,j)=a(n-2;i)+a(n-2;i+1)+\cdots+a(n-2;j-2)$.

Hence, using (\ref{eq1}) and (\ref{eq2}), we get the following theorem.

\begin{prop}
Let $s_n=|S_n(123,231)|$. For all $n\geq 3$,
$$s_n=s_{n-1}+s_n(1)+s_n(2)+\cdots+s_n(n-1),$$
where for all $1\leq i\leq n$,
$$s_n(i)=\sum_{j=1}^{i-1} s_{n-1}(j)+\sum_{j=i}^{n-2}(n-1-j)s_{n-2}(j),$$
and $s_3(1)=1$, $s_3(2)=1$, $s_3(3)=2$.
\end{prop}

Using this theorem, we get quickly the first values of the
sequence $|S_n(123,231)|$ for $n=0,1,2,\dots,10$:

\begin{center}
\begin{tabular}{|c|c|c|c|c|c|c|c|c|c|c|c|}
\hline
             $n$ & 0 &  1 &  2 &  3 &  4 &  5 &   6 &  7  &   8  &    9  &  10 \\
\hline
$|S_n(123,231)|$ & 1 &  1 &  2 &  4 & 11 & 39 & 161 & 784 & 4368 & 27260 & 189540 \\
\hline
\end{tabular}
\end{center}
\subsection{Avoidance of patterns 132 and 213 simultaneously}
We consider avoidance of the patterns 132 and 213 simultaneously.

Let $b(n;i_1,i_2,\dots,i_m)$ denote the number of permutations
$\pi\in S_n(132,213)$ such that $\pi_1\pi_2\dots\pi_m=i_1i_2\dots
i_m$ and let $b(n)=|S_n(132,213)|$.

Suppose $\pi\in S_n(132,213)$ is such that $\pi_1=i$ and
$\pi_2=j$. If $i>j$ then, since $\pi$ avoids $213$, we get
$\pi_3\leq i-1$. Thus
\begin{equation}
b(n;i,j)=\sum_{k=1,\ k\neq j}^{i-1} b(n-1;j,k).
\label{eq3}
\end{equation}

If $i<j$ then, since $\pi$ avoids $132$, we get $\pi_3\leq i-1$ or
$\pi_3\geq j+1$. Thus
\begin{equation}
b(n;i,j)=\sum_{k=1}^{i-1} b(n-1;j-1,k)+\sum_{k=j}^{n-1} b(n-1;j-1,k).
\label{eq4}
\end{equation}

Using (\ref{eq3}) and (\ref{eq4}), we get the following theorem.

\begin{prop}
Let $s_n=|S_n(132,213)|$. Then $s_n=\sum_{i,j=1}^n s(n;i,j)$ with

$s(n;i,i)=0$ for all $n,i\geq 1$;

$s(n;i,j)=\sum_{k=1}^{i-1} s(n-1;j,k)$ if $i>j$;

$s(n;i,j)=\sum_{k=1}^{i-1}s(n-1;j-1,k)+\sum_{k=j}^{n-1}s(n-1;j-1,k)$ if $i<j$;

and  $s(2;1,2)=s(2;2,1)=1$, $s(2;1,1)=s(2;1,1)=0$.
\end{prop}

Using this theorem, we get

\begin{center}
\begin{tabular}{|c|c|c|c|c|c|c|c|c|c|c|c|}
\hline
            $n$ & 0 &  1 &  2 &  3 &  4 &  5 &   6 &  7  &   8  &    9  &  10 \\
\hline
$|S_n(132,213)|$ & 1 &  1 &  2 &  4 & 11 & 37 & 149 & 705 & 3814 & 23199 & 156940 \\
\hline
\end{tabular}
\end{center}
\subsection{Avoidance of the patterns 213 and 231 simultaneously}
We now consider avoidance of the patterns 213 and 231
simultaneously. This case is equivalent to avoidance of the
patterns 132 and 312 by applying the reverse operation.

Let $c(n;i_1,i_2,\dots,i_m)$ denote the number of permutations
$\pi\in S_n(132,312)$ such that $\pi_1\pi_2\dots\pi_m=i_1i_2\dots
i_m$ and let $c(n)=|S_n(132,312)|$. We proceed as in the previous
case. For $n\geq i>j\geq 1$, we have
\begin{equation}
c(n;i,j)=\sum_{k=1}^{j-1}c(n-1;j,k)+\sum_{k=i}^{n-1}c(n-1;j,k).
\label{eq5}
\end{equation}

For $1\leq i<j\leq n$, we have
\begin{equation}
c(n;i,j)=\sum_{k=1}^{i-1}c(n-1;j-1,k)+\sum_{k=j}^{n-1}c(n-1;j-1,k).
\label{eq6}
\end{equation}

Using (\ref{eq5}) and (\ref{eq6}), we get the following theorem.

\begin{prop}
Let $s_n=|S_n(132,312)|$. Then $s_n=\sum_{i,j=1}^n s(n;i,j)$ with

$s(n;i,i)=0$ for all $n,i\geq 1$;

$s(n;i,j)=\sum_{k=1}^{j-1}s(n-1;j,k)+\sum_{k=i}^{n-1}s(n-1;j,k)$ if $i>j$;

$s(n;i,j)=\sum_{k=1}^{i-1}s(n-1;j-1,k)+\sum_{k=j}^{n-1}s(n-1;j-1,k)$ if $i<j$;

and $s(2;1,2)=s(2;2,1)=1$, $s(2;1,1)=s(2;1,1)=0$.
\end{prop}

Using this theorem, we get

\begin{center}
\begin{tabular}{|c|c|c|c|c|c|c|c|c|c|c|c|}
\hline
            $n$ & 0 &  1 &  2 &  3 &  4 &  5 &   6 &  7  &   8  &    9  &  10 \\
\hline
$|S_n(132,312)|$ & 1 &  1 &  2 &  4 & 10 & 30 & 108 & 454 & 2186 & 11840 & 71254 \\
\hline
\end{tabular}
\end{center}
\section{Avoiding a pattern x-y-z and beginning or ending with certain patterns}

Recall the definitions of $G_{q}^{p}(x)$, $N_{q}^{p}(n)$, $C(x)$ and $C_n$ in Section~\ref{Preliminaries}.

\begin{prop}\label{prop1} We have
    $$G_{1\mn3\mn2}^{12\ldots k}(x)=x^kC^2(x).$$
\end{prop}

\begin{proof}
Suppose $\pi={\pi}^{\prime}n{\pi}^{\prime\prime}\in
S_n(1\mn3\mn2)$ is such that $\pi_1 < \pi_2 < \cdots < \pi_k$ and
$\pi_j=n$. It is easy to see that $\pi$ avoids $1\mn3\mn2$ if and
only if ${\pi}^{\prime}$ is a $1\mn3\mn2$-avoiding permutation on
the letters $n-j+1,n-j+2,\ldots,n$, and ${\pi}^{\prime\prime}\in
S_{n-j}(1\mn3\mn2)$. If we now consider two cases, namely $j=k$
and $j\geq k+1$, we get
$$G_{1\mn3\mn2}^{12\ldots k}(x)=x^kC(x)+xG_{1\mn3\mn2}^{12\ldots k}(x)C(x).$$
Thus, $G_{1\mn3\mn2}^{12\ldots k}(x)=x^kC(x)/(1-xC(x))$ and,
using~(\ref{catalan}), we get the desired result.
\end{proof}

\begin{prop}\label{prop2} We have
$$G_{1\mn3\mn2}^{k(k-1)\ldots 1}(x)=x^kC^{k+1}(x).$$
\end{prop}
\begin{proof}
Suppose $\pi={\pi}^{\prime}n{\pi}^{\prime\prime}\in
S_n(1\mn3\mn2)$ is such that $\pi_1 > \pi_2 > \dots > \pi_k$ and
$\pi_j=n$. It is easy to see that $\pi$ avoids $1\mn3\mn2$ if and
only if ${\pi}^{\prime}$ is a $1\mn3\mn2$-avoiding permutation on
the letters $n-j+1, n-j+2, \ldots, n$, and
${\pi}^{\prime\prime}\in S_{n-j}(1\mn3\mn2)$. If we consider
separately the cases $j=1$ and $j\geq 2$, we get
$$G_{1\mn3\mn2}^{k(k-1)\ldots 1}(x)=xG_{1\mn3\mn2}^{(k-1)(k-2)\ldots 1}(x)+xG_{1\mn3\mn2}^{k(k-1)\ldots 1}(x)C(x).$$
Hence,
$$G_{1\mn3\mn2}^{k(k-1)\ldots 1}(x)=xG_{1\mn3\mn2}^{(k-1)(k-2)\ldots 1}(x)/(1-xC(x))$$
and, using~(\ref{catalan}), we get $G_{1\mn3\mn2}^{k(k-1)\ldots
1}(x)=xC(x)G_{1\mn3\mn2}^{(k-1)(k-2)\ldots 1}(x)$. By induction on
$k$, using the fact that $G_{1\mn3\mn2}^{1}(x)=C(x)-1=xC^2(x)$, we
get the desired result.
\end{proof}

\begin{prop}\label{prop3} We have
$$G_{2\mn1\mn3}^{12\ldots k}(x)=x^kC^{k+1}(x).$$
\end{prop}

\begin{proof}
One can use the same considerations as we have in the proof of
Proposition~\ref{prop2}, by considering a permutation
$\pi={\pi}^{\prime}1{\pi}^{\prime\prime}\in S_n(2\mn1\mn3)$ such
that $\pi_1 < \pi_2 < \dots < \pi_k$ and $\pi_j=1$.
\end{proof}

\begin{prop}\label{prop4} We have
$$G_{2\mn1\mn3}^{k(k-1) \ldots 1}(x)=x^kC^2(x).$$
\end{prop}
\begin{proof}
One can use the same considerations as we have in the proof of
Proposition~\ref{prop1}, by considering a permutation
$\pi={\pi}^{\prime}1{\pi}^{\prime\prime}\in S_n(2\mn1\mn3)$ such
that $\pi_1 > \pi_2 > \dots > \pi_k$ and $\pi_j=1$.
\end{proof}

Let $s_n(i_1,\dots,i_m)$ denote the number of permutations $\pi\in
S_n(1\mn2\mn3)$ such that $\pi_1\pi_2\ldots \pi_m = i_1i_2\ldots
i_m$. It is easy to see that
\begin{equation}
s_n(n)=s_n(n-1)=C_{n-1},
\label{eq7}
\end{equation}
and
\begin{equation}
s_n(t)=s_n(t,n) + \sum_{j=1}^{t-1} s_n(t,j) = s_{n-1}(t)+\sum_{j=1}^{t-1} s_{n-1}(j).
\label{eq8}
\end{equation}

Now, (\ref{eq7}) and (\ref{eq8}) with induction on $t$ give
\begin{equation}
s_n(n-t)=\sum_{j=0}^t (-1)^j\binom{t-j}{j}C_{n-j}
\label{eq9}
\end{equation}

Let us prove the following proposition.

\begin{prop}\label{prop5} We have
$$G_{1\mn2\mn3}^{12 \ldots k}(x)=\left\{
\begin{array}{ll}
0, & \mbox{if $k\geq 3$,} \\
x^2C^2(x), & \mbox{if $k=2$,}\\
xC^2(x), & \mbox{if $k=1$ .}
\end{array} \right.$$
\end{prop}
\begin{proof}
For $k\geq 3$, the statement is obviously true. If $k=1$ then
$G_{1\mn2\mn3}^{1}(x)=C(x)-1=xC^2(x)$.

Suppose now that $k=2$. From the definitions, for all $n\geq 2$, we have
$$N_{1\mn2\mn3}^{12}(n)=\sum_{i=1}^{n-1}\sum_{j=i+1}^n s_n(i,j).$$
In this formula, $j$ can only be equal to $n$, since otherwise we
have an occurrence of the pattern $1\mn2\mn3$. Using this fact
with (\ref{eq7}) and (\ref{eq8}), we get for $n\geq 2$,
$$N_{1\mn2\mn3}^{12}(n)=\sum_{i=1}^{n-1} s_n(i,n)=\sum_{i=1}^{n-1} s_{n-1}(i)=C_{n-1}.$$
Hence, $G_{1\mn2\mn3}^{12}(x)=x(C(x)-1)=x^2C^2(x)$.
\end{proof}

\begin{prop}\label{prop6} We have
$$N_{1\mn2\mn3}^{k(k-1) \ldots 1}(n)=
\sum_{t=1}^{n+1-k}\binom{n-t}{k-1}\sum_{j=0}^{n-t}(-1)^j\binom{n-t-j}{j}C_{n-t-j}.$$
\end{prop}
\begin{proof}
From the definitions, we have
$$N_{1\mn2\mn3}^{k(k-1)\ldots 1}(n)=\sum_{i_1=k}^n\sum_{i_2=1}^{i_1-1}\cdots\sum_{i_k=1}^{i_{k-1}-1}s_n(i_1,\dots,i_k)=\sum_{t=1}^{n+1-k} \binom{n-t}{k-1} s_n(t).$$
Using~(\ref{eq9}), we get
$$N_{1\mn2\mn3}^{k(k-1) \ldots 1}(n)=\sum_{t=1}^{n+1-k}\binom{n-t}{k-1}\sum_{j=0}^{n-t}(-1)^j\binom{n-t-j}{j}C_{n-t-j}.$$
\end{proof}

\section{Avoiding a pattern x-y-z, beginning and ending with certain patterns simultaneously}

Recall the definitions of $G_{q}^{p, r}(x)$ and $N_{q}^{p, r}(n)$
in Section~\ref{Preliminaries}.

\begin{prop}\label{prop10} We have

{\rm(i)} $G_{1\mn3\mn2}^{12\ldots k, 12\ldots \ell}(x) = x^{k+\ell
-1}C^{\ell+1}(x)+\frac{x^m-x^{k+\ell -1}}{1-x}$.

{\rm(ii)} $G_{1\mn3\mn2}^{12\ldots k, \ell(\ell-1) \ldots 1}(x) =
x^{k+\ell-1}C^2(x)$.

{\rm(iii)} $G_{1\mn3\mn2}^{k(k-1)\ldots 1, \ell(\ell-1)\ldots
1}(x) = x^{k+\ell-1}C^{k+1}(x)+\frac{x^m-x^{k+\ell -1}}{1-x}$,
where $m=max(k,\ell)$.

{\rm(iv)} the generating function
$G_{1\mn3\mn2}(x,y,z)=\sum_{k,\ell\geq0}G_{1\mn3\mn2}^{k(k-1)\ldots1,
12\ldots\ell}(x)y^kz^{\ell}$ for the sequence\\
$\{G_{1\mn3\mn2}^{k(k-1)\ldots1,12\ldots\ell}(x)\}_{k,\ell\geq0}$
{\rm(}where $k$ and $\ell$ go through all natural numbers{\rm)} is
$$\frac{1}{1-x(y+z)}\left( x(y+z+yz) +\frac{C(x)-1}{(1-xyC(x))(1-xzC(x))}\right).$$
\end{prop}
\begin{proof}
\

{\bf Beginning with $12\ldots k$ and ending with
$\ell(\ell-1)\ldots 1$:} Suppose
$\pi={\pi}^{\prime}n{\pi}^{\prime\prime}\in S_n(1\mn3\mn2)$ is
such that $\pi_1 < \pi_2 < \cdots < \pi_k$, $\pi_n < \pi_{n-1} <
\cdots < \pi_{n-\ell+1}$ and $\pi_j=n$. It is easy to see that
$\pi$ avoids $1\mn3\mn2$ if and only if ${\pi}^{\prime}$ is a
$1\mn3\mn2$-avoiding permutation on the letters
$n-j+1,n-j+2,\ldots,n$, and ${\pi}^{\prime\prime}\in
S_{n-j}(1\mn3\mn2)$. We now consider three cases, namely $j=k$,
$k+1\leq j \leq n-\ell$ and $j=n-\ell+1$. In terms of generating
functions, we have
$$
G_{1\mn3\mn2}^{12\ldots
k,\ell(\ell-1)\ldots1}(x)=x^kG_{2\mn1\mn3}^{\ell(\ell-1)\ldots
1}(x) + xG_{1\mn3\mn2}^{12\ldots
k}(x)G_{2\mn1\mn3}^{\ell(\ell-1)\ldots
1}(x)+x^{\ell}G_{1\mn3\mn2}^{12\ldots k}(x) + x^{k+\ell-1},$$
where we observed that to avoid $1\mn3\mn2$ and end with
$\ell(\ell-1)\ldots 1$ is the same as to avoid $2\mn1\mn3$ and
begin with $\ell(\ell-1)\ldots 1$ by applying the reverse and
complement operations. Also, we added the term $x^{k+\ell-1}$,
since when $j=k=n-\ell+1$, we have one ``good''
$(k+\ell-1)$-permutation, which is not counted by our three cases.

From Propositions~\ref{prop1} and~\ref{prop4}, we have that
$$G_{1\mn3\mn2}^{12\ldots k}(x)=x^kC^2(x)\mbox{ and }G_{2\mn1\mn3}^{\ell(\ell-1)\ldots1}(x)=x^{\ell}C^2(x).$$
Thus, using the fact that $xC^2(x)=C(x)-1$, we get
$$\begin{array}{l}
G_{1\mn3\mn2}^{12\ldots k,\ell(\ell-1)\ldots1}(x)=x^{k+\ell}C^2(x)(2+xC^2(x))+x^{k+\ell -1}\\
\qquad\qquad\qquad\qquad=x^{k+\ell-1}(C(x)-1)(C(x)+1)+x^{k+\ell-1}=x^{k+\ell-1}C^2(x).
\end{array}$$

{\bf Beginning with $12\ldots k$ and ending with $12\ldots \ell$:}
Suppose $\pi={\pi}^{\prime}n{\pi}^{\prime\prime}\in
S_n(1\mn3\mn2)$ is such that $\pi_1 < \pi_2 < \cdots < \pi_k$,
$\pi_n > \pi_{n-1} > \cdots > \pi_{n-\ell+1}$ and $\pi_j=n$. As
above, $\pi$ avoids $1\mn3\mn2$ if and only if ${\pi}^{\prime}$ is
a $1\mn3\mn2$-avoiding permutation on the letters
$n-j+1,n-j+2,\ldots,n$, and ${\pi}^{\prime\prime}\in
S_{n-j}(1\mn3\mn2)$. We consider the cases $j=k$, $k+1\leq j \leq
n-\ell$ and $j=n$. In terms of generating functions, the first
approximation for the function $G_{1\mn3\mn2}^{12\ldots k,12\ldots
\ell}(x)$ is
$$G_{1\mn3\mn2}^{12\ldots k,12\ldots \ell}(x) \approx x^kG_{2\mn1\mn3}^{12\ldots \ell}(x) + xG_{1\mn3\mn2}^{12\ldots k}(x)G_{2\mn1\mn3}^{12\ldots \ell}(x)+xG_{1\mn3\mn2}^{12\ldots k,12\ldots (\ell-1)}(x),$$
where we observed that to avoid $1\mn3\mn2$ and end with $12\ldots
\ell$ is the same as to avoid $2\mn1\mn3$ and begin with $12\ldots
\ell$ by applying the reverse and complement operations. We use
the sign ``$\approx$'' because there are some ``good''
permutations, which are not counted by our considerations. We
discuss them below.

From Propositions~\ref{prop1} and~\ref{prop3}, we have that
$G_{1\mn3\mn2}^{12\ldots k}(x)=x^kC^2(x)$ and
$G_{2\mn1\mn3}^{12\ldots \ell}(x)=x^{\ell}C^{\ell + 1}(x)$. Thus,
using the fact that $xC^2(x)=C(x)-1$ and $G_{1\mn3\mn2}^{12\ldots
k,1}(x) = G_{1\mn3\mn2}^{12\ldots k}(x) = x^kC^2(x)$
(Proposition~\ref{prop1}), we get
$$\begin{array}{l}
G_{1\mn3\mn2}^{12\ldots k,12\ldots \ell}(x)\\ \\
\qquad\approx x^{k+\ell}C^{\ell + 1}(x) + x^{k+\ell+1}C^{\ell +
3}(x)+xG_{1\mn3\mn2}^{12\ldots k,12\ldots (\ell-1)}(x)\\ \\
\qquad=x^{k+\ell}C^{\ell + 2}(x) + xG_{1\mn3\mn2}^{12\ldots
k,12\ldots (\ell-1)}(x)\\ \\
\qquad=x^{k+\ell}C^{\ell + 2}(x) + x^{k+\ell}C^{\ell + 1}(x) +
x^2G_{1\mn3\mn2}^{12\ldots k,12\ldots (\ell-2)}(x)\\ \\
\qquad=\cdots = x^{k+\ell}C^4(x)(C^{\ell-2}(x) +C^{\ell-3}(x)
+\cdots + 1)+x^{k+\ell-1}C^2(x)\\ \\
\qquad=x^{k+\ell-1}(C(x)-1)C^2(x)\frac{1-C^{\ell -
1}(x)}{1-C(x)}+x^{k+\ell-1}C^2(x)=x^{k+\ell-1}C^{\ell + 1}(x).
\end{array}$$
To complete the proof of this case, we observe that in our
considerations above, we do not count increasing permutations of
length $m=max(k,\ell)$, $m+1, \ldots, k+\ell-2$, which satisfy all
our restrictions. We did not count them because the $k$-beginning
and $\ell$-ending in these permutations overlap in more than one
letter. So, to get the desired result, we need to add the term
$$x^m + x^{m+1} + \ldots + x^{k+\ell-2}=(x^m-x^{k+\ell-1})/(1-x)$$
to the approximate value of $G_{1\mn3\mn2}^{12\ldots k,12\ldots
\ell}(x)$. For example, expanding $G_{1\mn3\mn2}^{12,123}(x)$, we
have, in particular, that there are 2002 10-permutations that
avoid $1\mn3\mn2$, begin with the pattern $12$ and end with the
pattern $123$.

{\bf Beginning with $k(k-1)\ldots 1$ and ending with
$\ell(\ell-1)\ldots 1$:} If $\ell = 1$ then, by
Proposition~\ref{prop2}, $G_{1\mn3\mn2}^{k(k-1)\ldots
1,1}(x)=x^kC^{k+1}(x)$. Suppose $\ell \geq 2$, and
$\pi={\pi}^{\prime}1{\pi}^{\prime\prime}\in S_n(1\mn3\mn2)$ is
such that $\pi_1 > \pi_2 > \cdots > \pi_k$, $\pi_n < \pi_{n-1} <
\cdots < \pi_{n-\ell+1}$ and $\pi_j=1$. Obviously,
${\pi}^{\prime\prime}$ is the empty word, since otherwise we have
an occurrence of the pattern $1\mn3\mn2$ starting from the letter
$1$. Thus, the first approximation for the function
$G_{1\mn3\mn2}^{k(k-1)\ldots 1,\ell(\ell-1)\ldots 1}$ is
$$G_{1\mn3\mn2}^{k(k-1)\ldots 1,\ell(\ell-1)\ldots 1}(x) \approx xG_{1\mn3\mn2}^{k(k-1)\ldots 1,(\ell-1)(\ell-2)\ldots 1}(x)=\cdots =x^{k+\ell-1}C^{k+1}(x).$$
Like in the previous case, we did not count decreasing
permutations of length $m=max(k,\ell)$, $m+1, \ldots, k+\ell-2$,
which satisfy all our restrictions. Thus, to get the desired
result, we add the term $(x^m-x^{k+\ell-1})/(1-x)$ to the
approximate value of $G_{1\mn3\mn2}^{k(k-1)\ldots
1,\ell(\ell-1)\ldots 1}(x)$.

{\bf Beginning with $k(k-1)\ldots 1$ and ending with $12\ldots
\ell$:} Suppose $\pi={\pi}^{\prime}n{\pi}^{\prime\prime}\in
S_n(1\mn3\mn2)$. Any letter of ${\pi}^{\prime}$ is greater than
any letter of ${\pi}^{\prime\prime}$, since otherwise we have an
occurrence of the pattern $1\mn3\mn2$ in $\pi$ containing the
letter $n$ which is forbidden. Also, ${\pi}^{\prime}$ and
${\pi}^{\prime\prime}$ avoid $1\mn3\mn2$. If $\pi$ begins with
$k(k-1)\ldots 1$, ends with $12\ldots \ell$ and ${\pi}^{\prime}$
and ${\pi}^{\prime\prime}$ are not empty, then ${\pi}^{\prime}$
must begin with $k(k-1)\ldots 1$ and ${\pi}^{\prime\prime}$ must
end with $12\ldots \ell$. If ${\pi}^{\prime}$ is empty then
${\pi}^{\prime\prime}$ must begin with $(k-1)(k-2)\ldots 1$ and
end with $12\ldots \ell$. If ${\pi}^{\prime\prime}$ is empty then
${\pi}^{\prime}$ must begin with $k(k-1)\ldots 1$ and end with
$12\ldots (\ell-1)$. In terms of generating functions, the
discussion above leads to the following:
$$G_{1\mn3\mn2}^{k(k-1) \ldots 1, 12\ldots \ell}(x)\approx
xG_{1\mn3\mn2}^{k(k-1) \ldots 1}(x)G_{2\mn1\mn3}^{12\ldots
\ell}(x) + xG_{1\mn3\mn2}^{(k-1)\ldots 1, 12\ldots \ell}(x) +
xG_{1\mn3\mn2}^{k(k-1) \ldots 1, 12\ldots (\ell-1)}(x),$$ where we
observed that to avoid $1\mn3\mn2$ and end with $12\ldots \ell$ is
the same as to avoid $2\mn1\mn3$ and begin with $12\ldots \ell$.
However, to put the sign ``$=$'' instead of ``$\approx$'', we have
to correct the right-hand side of the recurrence relation by
observing that when either $k=1$ and $\ell=0$, or $k=0$ and $\ell
= 1$, or $k=1$ and $\ell=1$, the formula do not count the
permutation $\pi=1$ which satisfies all the conditions needed.
Thus, if we make correction of the right-hand side, then multiply
both parts of the obtained equality by $x^ky^{\ell}$ and sum over
all natural $k$ and $\ell$ we get (recall the definition of
$G_{1\mn3\mn2}(x,y,z)$ in the statement of the theorem):
$$G_{1\mn3\mn2}(x,y,z)=x\displaystyle\sum_{k,\ell \geq 0}G_{1\mn3\mn2}^{k(k-1) \ldots
1}(x)G_{2\mn1\mn3}^{12\ldots \ell}(x)y^kz^{\ell} +
x(y+z)G_{1\mn3\mn2}(x,y,z) + x(y+z+yz).$$ From Propositions 5 and
6, $G_{1\mn3\mn2}^{k(k-1) \ldots 1}(x)G_{2\mn1\mn3}^{12\ldots
\ell}(x) = x^{k+\ell}C^{k+\ell+2}(x)$, and thus
$$\begin{array}{l}
G_{1\mn3\mn2}(x,y,z)=\frac{1}{1-x(y+z)}\left(x(y+z+yz)+\displaystyle\sum_{k,\ell
\geq 0}x^{k+\ell}C^{k+\ell+2}(x)y^kz^{\ell} \right)\\
\qquad\qquad\qquad\,=\frac{1}{1-x(y+z)}\left(
x(y+z+yz)+zC^2(z)\displaystyle\sum_{k \geq
0}(xyC(x))^k\displaystyle\sum_{\ell \geq 0}(xzC(x))^{\ell}
\right)\\
\qquad\qquad\qquad\,=\frac{1}{1-x(y+z)}\left( x(y+z+yz) +
\frac{C(x)-1}{(1-xyC(x))(1-xzC(x))}\right), \end{array}$$ where we
used that $xC^2(x)=C(x)-1$.
\end{proof}

\begin{prop}\label{prop11} We have

{\rm(i)} $G_{2\mn1\mn3}^{12\ldots k, 12\ldots \ell}(x) = x^{k+\ell
-1}C^{k+1}(x)+\frac{x^m-x^{k+\ell -1}}{1-x}$.

{\rm(ii)} $G_{2-1-3}^{k(k-1)\ldots 1, 12\ldots \ell}(x) =
x^{k+\ell -1}C^2(x)$.

{\rm(iii)} $G_{2\mn1\mn3}^{k(k-1)\ldots 1, \ell(\ell-1)\ldots
1}(x) = x^{k+\ell -1}C^{\ell+1}(x)+\frac{x^m-x^{k+\ell -1}}{1-x}$,
where $m=max(k,\ell)$.

{\rm(iv)} the generating function
$G_{2\mn1\mn3}(x,y,z)=\sum_{k,\ell \geq 0} G_{2\mn1\mn3}^{12\ldots
k, \ell(\ell-1) \ldots 1}(x)y^kz^{\ell}$ for the sequence\\
$\{G_{2\mn1\mn3}^{12\ldots k, \ell(\ell-1) \ldots
1}(x)\}_{k,\ell\geq0}$ {\rm(}where $k$ and $\ell$ go through all
natural numbers{\rm)} is
$$\frac{1}{1-x(y+z)}\left( x(y+z+yz) + \frac{C(x)-1}{(1-xyC(x))(1-xzC(x))}\right).$$
\end{prop}
\begin{proof}
We apply the reverse and complement operations and then use the
results of Proposition~\ref{prop10}. For example, to avoid
$2\mn1\mn3$, begin with $12\ldots k$ and end with $12\ldots \ell$
is the same as to avoid $1\mn3\mn2$, begin with $12\ldots \ell$
and end with $12\ldots k$.
\end{proof}

Let $h_n^{k,\ell}(t;s)$ denote the number of $1\mn2\mn3$-avoiding
$n$-permutations such that $\pi_k=t$, $\pi_{n-\ell+1}=s$,
$\pi_1>\pi_2>\cdots>\pi_k$, and
$\pi_{n-\ell+1}>\pi_{n-\ell+2}>\cdots>\pi_n$. Also, we define
$g_n(i_1,i_2,\ldots,i_m;b)$ to be the number of
$1\mn2\mn3$-avoiding $n$-permutations such that
$\pi_1\pi_2\cdots\pi_m=i_1i_2\ldots i_m$ and $\pi_n=b$. We need
the following two lemmas to prove Proposition~\ref{prop12}.

\begin{lemma}\label{aux0} For all $n\geq 2$,
$$g_n(a;b)=\left\{
\begin{array}{ll}
0, &             2\leq a+1<b\leq n,\\
\binom{n-2}{a-1}, & 1\leq a\leq n-1,\\
\sum\limits_{j=0}^{n-a}(-1)^j\binom{n-a-j}{j}\left(\sum\limits_{i=0}^{b-1}(-1)^i\binom{b-1-i}{i}C_{n-2-j-i}
\right), & 1\leq b<a\leq n.
\end{array}\right.$$
\end{lemma}

\begin{proof}
By definitions we have

(1) $g_n(a;b)=0$ for all $2\leq a+1<b\leq n$;

(2)
$g_n(a;a+1)=g_n(a,1;a+1)+\ldots+g_n(a,a-1;a+1)+g_n(a,a+2;a+1)+\ldots+g_n(a,n-1;a+1)+g_n(a,n;a+1)$.
Using the fact that no there exists a permutation $\pi\in
S_n(1\mn2\mn3)$ such that $\pi_1=a$, $\pi_2\leq a-2$, and
$\pi_n=a+1$ we get
$$g_n(a;a+1)=g_n(a,a-1;a+1)+g_n(a,a+2;a+1)+\ldots+g_n(a,n;a+1).$$
Using the fact that no there exists a permutation $\pi\in
S_n(1\mn2\mn3)$ such that $\pi_1=a$ and $a\leq\pi_2\leq n-1$ we
get $g_n(a;a+1)=g_n(a,a-1;a+1)+g_n(a,n;a+1)$. On the other hand,
it is easy to see that $g_n(a,a-1;a+1)=g_{n-1}(a-1;a)$ and
$g_n(a,n;a+1)=g_{n-1}(a;a+1)$. Hence,
                  $$g_n(a;a+1)=g_{n-1}(a-1;a)+g_{n-1}(a;a+1).$$
Using induction we get that $g_n(a;a+1)=\binom{n-2}{a-1}$ for all
$n\geq 2$ and $1\leq a\leq n-1$.

(3) Similarly as (2) we have for all $a>b$,
     $$g_n(a;b)=g_{n-1}(b-1;b)+g_{n-1}(b+1;b)+g_{n-1}(b+2;b)+\cdots+g_{n-1}(a;b).$$
Using Equation~(\ref{eq9}) we get
     $$g_n(a;1)=g_n(a;2)=s_{n-1}(a-1)=\sum_{j=0}^{n-a} (-1)^j\binom{n-a-j}{j} C_{n-2-j}.$$
Using induction on $b$, we get
     $$g_n(a;b)=\sum_{j=0}^{n-a}(-1)^j\binom{n-a-j}{j}\left(\sum_{i=0}^{b-1}(-1)^{i}\binom{b-1-i}{i}C_{n-2-j-i} \right).$$
\end{proof}

\begin{lemma}\label{aux} We have
$$h_n^{k,\ell}(t;s)=\left\{
\begin{array}{ll}
  \binom{n-t}{k-1}\binom{s-1}{\ell-1}g_{t-(\ell-1);s-(\ell-1)}(n+2-k-\ell),  & \mbox{if } 1\leq s<t\leq n;\\
  h_n^{k,\ell}(t+1;t),    & \mbox{if } s=t+1;\\
  h_{n-1}^{k,\ell-1}(t;s-1)+h_{n-1}^{k-1;\ell}(t;s-1),    & \mbox{if } 2\leq t+1<s\leq n. \\
\end{array}
\right.$$
\end{lemma}
\begin{proof}

(1) Let $n\geq t>s\geq 1$; so by definitions we get
  $$h_n^{k,\ell}(t;s)=\binom{n-t}{k-1}\binom{s-1}{\ell-1}g_{t-(\ell-1);s-(\ell-1)}(n-(k-1)-(\ell-1)).$$

(2) Let $s=t+1$; so it is easy to see $h_n^{k,\ell}(t;t+1)=h_n^{k,\ell}(t+1;t)$;

(3) Let $2\leq t+1<s\leq n$. Let $\pi$ any permutations in
$S_n(1\mn2\mn3)$ such that $\pi_k=t$ and $\pi_{n+1-\ell}=s$ where
$\pi_1>\cdots>\pi_k$ and $\pi_{n+1-\ell}>\cdots>\pi_n$; so there
two possibilities either $\pi_{n+2-\ell}=s-1$ or $\pi_j=s-1$ where
$j\leq k-1$. In this first case we get that there exist
$h_{n-1}^{k,\ell-1}(t;s-1)$ permutations, and in the second case
we have that there exist $h_{n-1}^{k-1;\ell}(t;s-1)$ permutations.
(we extend the number $h_n^{k,l}(a;b)$ as $0$ for any $\ell\leq 0$
or $k\leq 0$.
\end{proof}

We recall that the Kronecker delta $\delta_{n,k}$ is defined to be
\[ \delta_{n,k} = \left\{ \begin{array}{ll} 1, & \mbox{if $n=k$,}
\\ 0, & \mbox{else.}
\end{array}
\right. \]

\begin{prop}\label{prop12} We have

{\rm(i)} $G_{1\mn2\mn3}^{12\ldots k, 12\ldots \ell}(x)= \left\{
\begin{array}{ll} 0, & \mbox{if $k\geq 3$ or $\ell \geq 3$} \\
xC^2(x), & \mbox{if $k=1$ and $\ell=1$}
\end{array}
\right.$, $N_{1\mn2\mn3}^{12,12}(n)=\left\{ \begin{array}{ll} 0, &
\mbox{if $n=3$} \\ C_{n-2}, & \mbox{else}
\end{array}
\right.$, and\\
$N_{1\mn2\mn3}^{12,1}(n)=N_{1-2-3}^{1,12}(n)=C_{n-1}$.

{\rm(ii)} $N_{1\mn2\mn3}^{k(k-1)\ldots 1, 12\ldots\ell}(n)=\left\{
\begin{array}{ll} 0, & \mbox{if
$\ell \geq 3$,} \\
\sum\limits_{t=1}^{n-k}\binom{n-t-1}{k-1}\sum\limits_{j=0}^{n-t-1}(-1)^j\binom{n-t-j-1}{j}C_{n-t-j-1}
+ (k-1)\delta_{n,k+1}, & \mbox{if $\ell = 2$,} \\
\sum\limits_{t=1}^{n+1-k}\binom{n-t}{k-1}\sum\limits_{j=0}^{n-t}(-1)^j\binom{n-t-j}{j}C_{n-t-j},
& \mbox{if $\ell = 1$.}
\end{array}\right.$

{\rm(iii)} $N_{1\mn2\mn3}^{12\ldots k,
\ell(\ell-1) \ldots 1}(n)=\left\{ \begin{array}{ll} 0, & \mbox{if $k\geq 3$,} \\
\sum\limits_{t=1}^{n-\ell}\binom{n-t-1}{\ell-1}\sum\limits_{j=0}^{n-t-1}(-1)^j\binom{n-t-j-1}{j}C_{n-t-j-1}
+ (\ell-1)\delta_{n,\ell+1}, & \mbox{if $k=2$,} \\
\sum\limits_{t=1}^{n+1-\ell}\binom{n-t}{\ell-1}\sum\limits_{j=0}^{n-t}(-1)^j\binom{n-t-j}{j}C_{n-t-j},
& \mbox{if $k=1$.}
\end{array}\right.$

{\rm(iv)} $N_{1-2-3}^{k(k-1)\ldots 1, \ell(\ell-1)\ldots 1}(x) =
\sum_{t=1}^{n-k+1}\sum_{s=\ell}^{n} h_n^{k,\ell}(t;s)$,  where
$h_n^{k,\ell}(t;s)$ is given in Lemma~\ref{aux}.
\end{prop}
\begin{proof}
\

{\bf Beginning with $12\ldots k$ and ending with $12\ldots \ell$:}
If $k\geq 3$ or $\ell \geq 3$, the statement is obvious, since in
this case $12\ldots k$ or $12\ldots \ell$ does not avoid the
pattern $1\mn2\mn3$. If $k=1$ or $\ell = 1$, we get the  statement
from Proposition~\ref{prop5} (in the first of these cases we apply
the reverse and complement operations). Suppose now that $k=2$,
$\ell = 2$, and an $n$-permutation $\pi$ avoids $1\mn2\mn3$,
begins with the pattern $12$ and ends with the pattern $12$. The
letter $n$ must be next to the leftmost letter, since otherwise
two leftmost letters and $n$ form the pattern $1\mn2\mn3$. Also,
the letter $1$ must be next to the rightmost letter, since
otherwise $1$ and two rightmost letters form the pattern
$1\mn2\mn3$. It is easy to see now that there are $C_{n-2}$
possibilities to choose $\pi$, since we can take any
$1\mn2\mn3$-avoiding permutation on the letters $\{2, 3, \ldots,
n-1\}$ (there are $C_{n-2}$ such permutations), then let the
letters $n$ and $1$ be in the second and $(n-1)$-st positions
respectively. These considerations are fail only when $n=3$, since
in this case the second and $(n-1)$-st positions coincide.
However, in this case we obviously have no permutations with the
good properties.

{\bf Beginning with $k(k-1)\ldots 1$ and ending with $12\ldots
\ell$:} The statement is true for $\ell \geq 3$, since in this
case $12\ldots \ell$ does not avoid $1\mn2\mn3$. For the case
$\ell = 1$ we use Proposition~\ref{prop6}. Suppose now that $\ell
= 2$, and an $n$-permutation $\pi$ avoids $1\mn2\mn3$, begins with
the pattern $k(k-1)\ldots 1$ and ends with the pattern $12$. The
letter $1$ must be next to the rightmost letter, since otherwise
$1$ and two rightmost letters form the pattern $1\mn2\mn3$. So, to
form $\pi$ we can take any $(n-1)$-permutation on the letters $\{
2, 3, \ldots, n \}$ that avoid $1\mn2\mn3$ and begin with the
pattern $k(k-1)\ldots 1$ (the number of such permutations is given
by Proposition~\ref{prop6}), and then let the letter $1$ be in the
$(n-1)$-st position. Also, we observe that in the case $n=k+1$ we
have $k-1$ extra permutations, which are obtained from the
$(n-1)$-permutations having the $k-1$ leftmost letters in
decreasing order and two rightmost letters in increasing order.

{\bf Beginning with $12\ldots k$ and ending with
$\ell(\ell-1)\ldots 1$:} By the reverse and complement operations,
to avoid $1\mn2\mn3$, begin with the pattern $12\ldots k$ and end
with the pattern $\ell(\ell-1)\ldots 1$ is the same as to avoid
$1\mn2\mn3$, begin with the pattern $\ell(\ell-1)\ldots 1$ and end
with the pattern $12\ldots k$, so we can apply the results of the
previous case.

{\bf Beginning with $k(k-1)\ldots 1$ and ending with
$\ell(\ell-1)\ldots 1$:} The statement is straitforward  to prove.
\end{proof}

\section{Avoiding a pattern xyz, beginning and ending with certain patterns simultaneously}

Recall the definitions of $E_{q}^{p, r}(x)$ in Section~\ref{Preliminaries}.

\begin{prop} We have

{\rm(i)}
$E_{213}^{12\ldots k, 12\ldots \ell}(x) = \left\{ \begin{array}{ll} E_{132}^{12\ldots \ell}(x), & \mbox{if $k=1$} \\
                                                                    E_{213}^{12\ldots k}(x), & \mbox{if $\ell=1$} \\
\end{array}\right.$,
where $E_{132}^{12\ldots \ell}(x)$ and $E_{213}^{12\ldots k}(x)$
are given in~\cite[Theorem 6]{Kit3} and~\cite[Theorem 10]{Kit3}
respectively. For $k,\ell \geq 2$, $E_{213}^{12\ldots k, 12\ldots
\ell}(x)$ satisfies
$$\frac{\partial}{\partial x}E_{213}^{12 \ldots k, 12\ldots \ell}(x) = E_{213}^{12\ldots k, 12\ldots (\ell-1)}(x) + \left(E_{213}^{12\ldots k,12}(x) + \frac{x^{k-1}}{(k-1)!} \right)E_{132}^{12\ldots \ell}(x).$$

{\rm(ii)}
$E_{213}^{12\ldots k,\ell(\ell-1) \ldots 1}(x) = \left\{ \begin{array}{ll}  E_{132}^{\ell(\ell-1) \ldots 1}(x), & \mbox{if $k=1$} \\
                                                                            E_{213}^{12\ldots k}(x), & \mbox{if $\ell=1$} \\
\end{array}\right.$,
where $E_{132}^{\ell(\ell-1) \ldots 1}(x)$ and $E_{213}^{12\ldots
k}(x)$ are given in~\cite[Theorem 7]{Kit3} and~\cite[Theorem
10]{Kit3} respectively. For $k,\ell \geq 2$, $E_{213}^{12\ldots
k,\ell(\ell-1) \ldots 1}(x)$ satisfies
$$\begin{array}{l}
\frac{\partial}{\partial x}E_{213}^{12 \ldots k, \ell(\ell-1)
\ldots 1}(x)\\
\qquad\qquad= \frac{x^{\ell-1}}{(\ell-1)!}E_{213}^{12\ldots k}(x)
+ \left( E_{213}^{12\ldots k, 12}(x)+
\frac{x^{k-1}}{(k-1)!}\right) E_{132}^{\ell (\ell-1)\ldots 1}(x)+
\binom{k+\ell
-2}{k-1}\frac{x^{k+\ell-2}}{(k+\ell-2)!}.\end{array}$$

{\rm(iii)}
$E_{213}^{k(k-1)\ldots 1, 12\ldots \ell}(x) = \left\{ \begin{array}{ll}
E_{132}^{12\ldots \ell}(x), & \mbox{if $k=1$} \\
E_{213}^{k(k-1)\ldots 1}(x), & \mbox{if $\ell=1$} \\
\end{array}\right.$,
where $E_{132}^{12\ldots \ell}(x)$ and $E_{213}^{k(k-1)\ldots
1}(x)$ are given in~\cite[Theorem 6]{Kit3} and~\cite[Theorem
11]{Kit3} respectively. For $k,\ell \geq 2$,
$E_{213}^{k(k-1)\ldots 1, 12\ldots \ell}(x)$ satisfies
$$\frac{\partial}{\partial x}E_{213}^{k(k-1)\ldots 1, 12\ldots
\ell}(x)= E_{213}^{(k-1)\ldots 1, 12\ldots \ell}(x) +
E_{213}^{k(k-1)\ldots 1,12}(x)E_{132}^{12\ldots \ell}(x) +
E_{213}^{k(k-1)\ldots 1, 12\ldots (\ell-1)}(x).$$

{\rm(iv)} $E_{213}^{k(k-1)\ldots 1, \ell(\ell-1)\ldots 1}(x) =
\left\{ \begin{array}{ll}
E_{132}^{\ell(\ell-1)\ldots 1}(x), & \mbox{if $k=1$} \\
E_{213}^{k(k-1)\ldots 1}(x), & \mbox{if $\ell=1$} \\
\end{array}\right.$,
where $E_{132}^{\ell(\ell-1)\ldots 1}(x)$ and
$E_{213}^{k(k-1)\ldots 1}(x)$ are given in~\cite[Theorem 7]{Kit3}
and~\cite[Theorem 11]{Kit3} respectively. For $k, \ell \geq 2$,
$E_{213}^{k(k-1)\ldots 1, \ell(\ell-1)\ldots 1}(x)$ satisfies
$$\frac{\partial}{\partial x}E_{213}^{k(k-1)\ldots 1,
\ell(\ell-1)\ldots 1}(x)=E_{213}^{(k-1)\ldots 1,
\ell(\ell-1)\ldots 1}(x) + \left(E_{132}^{\ell(\ell-1)\ldots 1}(x)
+ \frac{x^{\ell -1}}{(\ell - 1)!}\right)E_{213}^{k(k-1)\ldots 1,
12}(x).$$
\end{prop}

\begin{proof}
\

{\bf Beginning with $12\ldots k$ and ending with
$\ell(\ell-1)\ldots 1$:} The statement is obviously true when
$k=1$ and $\ell=1$. Suppose now that $k\geq 2$, $\ell \geq 2$ and
an $(n+1)$-permutation $\pi$ avoids $213$, begins with the pattern
$12\ldots k$ and ends with the pattern $12\ldots \ell$. The letter
$(n+1)$ can only be in the position $k$, or in the position $i$,
where $(k+1)\leq i \leq n-\ell+1$, or in the position $n-\ell+2$.
In the first case, we choose the $(k-1)$ leftmost letters in
$\binom{n}{k-1}$ ways, rearrange them into the increasing order,
and observe, that the letters of $\pi$ to the right of $(n+1)$
must form an $(n-k+1)$-permutation, that avoids $213$ and ends
with the pattern $\ell(\ell-1)\ldots 1$ (the number of such
permutations, using the reverse and complement operation, is equal
to the number of $(n-k+1)$-permutations that avoid $132$ and begin
with the pattern $\ell(\ell-1)\ldots 1$). In the third case, we
choose the $(\ell-1)$ rightmost letters in $\binom{n}{\ell-1}$
ways, rearrange them into the decreasing order, and observe, that
the letters of $\pi$ to the right of $(n+1)$ must form an
$(n-\ell+1)$-permutation, that avoids $213$, begins with the
pattern $12\ldots k$, and ends with the pattern $12$ (if it ends
with the pattern $21$, the letter $(n+1)$ and two letters
immediately to the left of it form the pattern $213$). In the
second case, we choose the letters of $\pi$ to the left of $(n+1)$
in $\binom{n}{i-1}$ ways and observe, that these letters must form
a $(i-1)$-permutation that avoids $213$, begins with the pattern
$12\ldots k$ and ends with the pattern $12$. In the same time, the
letters to the right of $(n+1)$ must form an $(n-i+2)$-permutation
that avoids $213$ and ends with the pattern $\ell(\ell-1)\ldots
1$. Besides, we observe that if $n=k+\ell -2$, that is
$|\pi|=k+\ell -1$, and first $k$-letters of $\pi$ are rearranged
into the increasing order, whereas the last $\ell$ letters are
rearranged in the decreasing order, we have a number of extra
``good'' permutations. The number of such permutations is the
number of ways of choosing the first $(k-1)$ letters, that is
$\binom{k+\ell-2}{k-1}$. This discussion leads to the following:
$$\begin{array}{l}
N_{213}^{12\ldots k,\ell(\ell-1)\ldots
1}(n+1)=\binom{n}{k-1}N_{132}^{\ell(\ell-1)\ldots
1}(n-k+1)+\binom{n}{\ell-1}N_{213}^{12\ldots k}(n-\ell+1)\\
\qquad\qquad\qquad\quad\qquad\qquad+\displaystyle\sum_{i=0}^{n}\binom{n}{i}N_{213}^{12\ldots
k,12}(i)N_{132}^{\ell(\ell-1)\ldots
1}(n-i)+\binom{k+\ell-2}{k-1}{\delta}_{n,k+\ell-2}, \end{array}$$
where $\delta_{n,k+\ell-2}$ is the Kronecker delta. We get the
desirable result by multiplying both sides of the last equality by
$x^n/n!$ and summing over $n$.

{\bf Beginning with $12\ldots k$ and ending with $12\ldots \ell$:}
The statement is obviously true when $k=1$ and $\ell=1$. Suppose
now that $k\geq 2$, $\ell \geq 2$ and an $(n+1)$-permutation $\pi$
avoids $213$, begins with the pattern $12\ldots k$ and ends with
the pattern $12\ldots \ell$. The letter $(n+1)$ can only be in the
position $k$, or in the position $i$, where $(k+1)\leq i \leq
n-\ell$, or in the $(n+1)$-th position. In the last case, the
number of such permutations is obviously $N_{213}^{12\ldots k,
12\ldots \ell-1}(n)$. In the first case, we choose the $(k-1)$
leftmost letters in $\binom{n}{k-1}$ ways, rearrange them into
increasing order, and observe, that the letters of $\pi$ to the
right of $(n+1)$ must form an $(n-k+1)$-permutation, that avoids
$213$ and ends with the pattern $12\ldots \ell$ (the number of
such permutations, using the reverse and complement operation, is
equal to the number of $(n-k+1)$-permutations that avoid $132$ and
begin with the pattern $12\ldots \ell$). In the second case, we
choose the letters of $\pi$ to the left of $(n+1)$ in
$\binom{n}{i-1}$ ways and observe, that these letters must form a
$(i-1)$-permutation that avoids $213$, begins with the pattern
$12\ldots k$ and ends with the pattern $12$ (if it ends with the
pattern $21$, the letter $(n+1)$ and two letters immediately to
the left of it form the pattern $213$). In the same time, the
letters to the right of $(n+1)$ must form an $(n-i+2)$-permutation
that avoids $213$ and ends with the pattern $12\ldots \ell$. This
discussion leads to the following:
$$\begin{array}{l}
N_{213}^{12\ldots k,12\ldots \ell}(n+1)=N_{213}^{12\ldots k,
12\ldots \ell-1}(n)\\
\qquad\qquad\qquad\qquad\qquad+\displaystyle\sum_{i=0}^{n}\binom{n}{i}N_{213}^{12\ldots
k,12}(i)N_{132}^{12\ldots \ell}(n-i) +
\binom{n}{k-1}N_{132}^{12\ldots \ell}(n-k+1).\end{array}$$ We get
the desirable result by multiplying both sides of the last
equality by $x^n/n!$ and summing over $n$.

{\bf Beginning with $k(k-1)\ldots 1$ and ending with $12\ldots \ell$ or with $\ell(\ell-1)\ldots 1$:}
We proceed in the same way as we do under considering the previous case.
\end{proof}

\begin{prop}\label{prop13} We have

{\rm(i)}
$E_{132}^{12\ldots k, 12\ldots \ell}(x) =
\left\{ \begin{array}{ll} E_{213}^{12\ldots \ell}(x), & \mbox{if $k=1$} \\
E_{132}^{12\ldots k}(x), & \mbox{if $\ell=1$} \\
\end{array}\right.$,
where $E_{213}^{12\ldots \ell}(x)$ and $E_{132}^{12\ldots k}(x)$
are given in~\cite[Theorem 10]{Kit3} and~\cite[Theorem 6]{Kit3}
respectively. For $k, \ell \geq 2$, $E_{132}^{12\ldots k, 12\ldots
\ell}(x)$ satisfies
$$\frac{\partial}{\partial x}E_{132}^{12 \ldots k, 12\ldots \ell}(x) = E_{132}^{12\ldots k-1, 12\ldots \ell}(x) + \left(E_{132}^{12, 12\ldots \ell}(x) + \frac{x^{\ell-1}}{(\ell-1)!} \right)E_{132}^{12\ldots k}(x).$$

{\rm(ii)}
$E_{132}^{12\ldots k, \ell(\ell-1) \ldots 1}(x) =
\left\{ \begin{array}{ll} E_{132}^{12\ldots k}(x), & \mbox{if $\ell =1$} \\
E_{213}^{\ell (\ell -1)\ldots 1}(x), & \mbox{if $k=1$} \\
\end{array}\right.$,
where $E_{132}^{12\ldots k}(x)$ and $E_{213}^{\ell (\ell -1)
\ldots 1}(x)$ are given in~\cite[Theorem 6]{Kit3}
and~\cite[Theorem 11]{Kit3} respectively. For $k,\ell \geq 2$,
$E_{132}^{12\ldots k, \ell(\ell-1) \ldots 1}(x)$ satisfies
$$\frac{\partial}{\partial x}E_{132}^{12\ldots k, \ell(\ell-1)
\ldots 1}(x)= E_{132}^{12\ldots k, (\ell-1)\ldots 1}(x) +
E_{132}^{12, \ell(\ell-1)\ldots 1}(x)E_{132}^{12\ldots k}(x) +
E_{132}^{12 \ldots (k-1), \ell (\ell -1) \ldots 1}(x).$$

{\rm(iii)} $E_{132}^{k(k-1)\ldots 1, 12\ldots \ell}(x)=
\left\{ \begin{array}{ll} E_{213}^{12\ldots \ell}(x), & \mbox{if $k=1$} \\
E_{132}^{k(k-1)\ldots 1}(x), & \mbox{if $\ell=1$} \\
\end{array}\right.$,
where $E_{213}^{12\ldots \ell}(x)$ and $E_{132}^{k(k-1)\ldots
1}(x)$ are given in~\cite[Theorem 10]{Kit3} and~\cite[Theorem
7]{Kit3} respectively. For $k,\ell \geq 2$, $E_{132}^{k(k-1)\ldots
1, 12\ldots \ell}(x)$ satisfies
$$\begin{array}{l}
\frac{\partial}{\partial x}E_{213}^{k(k-1)\ldots 1, 12\ldots
\ell}(x) = \frac{x^{k-1}}{(k-1)!}E_{213}^{12\ldots \ell}(x)\\
\qquad\qquad\qquad\qquad\qquad+\left( E_{132}^{12,12\ldots
\ell}(x)+\frac{x^{\ell-1}}{(\ell-1)!}\right)E_{132}^{k(k-1)\ldots
1}(x)+ \binom{k+\ell -2}{\ell-1}\frac{x^{k+\ell-2}}{(k+\ell
-2)!}.\end{array}$$

{\rm(iv)}
$E_{132}^{k(k-1)\ldots 1, \ell(\ell-1)\ldots 1}(x) =
\left\{ \begin{array}{ll} E_{213}^{\ell(\ell-1)\ldots 1}(x), & \mbox{if $k=1$} \\
E_{132}^{k(k-1)\ldots 1}(x), & \mbox{if $\ell=1$} \\
\end{array}\right.$,
where $E_{132}^{k(k-1)\ldots 1}(x)$ and
$E_{213}^{\ell(\ell-1)\ldots 1}(x)$ are given in~\cite[Theorem
7]{Kit3} and~\cite[Theorem 11]{Kit3} respectively. For $k,\ell
\geq 2$, $E_{132}^{k(k-1)\ldots 1, \ell(\ell-1)\ldots 1}(x)$
satisfies
$$
\frac{\partial}{\partial x}E_{132}^{k(k-1)\ldots 1, \ell (\ell
-1)\ldots 1}(x)=E_{132}^{(\ell-1)\ldots 1, k(k-1)\ldots 1}(x) +
\left(E_{132}^{k(k-1)\ldots 1}(x) + \frac{x^{k -1}}{(k -
1)!}\right)E_{132}^{12, \ell(\ell-1)\ldots 1}(x).$$
\end{prop}

\begin{proof}
We apply the reverse and complement operations and then use the
results of Proposition~\ref{prop13}. For example, to avoid $213$,
begin with $12\ldots k$ and end with $12\ldots \ell$ is the same
as to avoid $132$, begin with $12\ldots \ell$ and end with
$12\ldots k$.
\end{proof}

\begin{prop} We have

{\rm(i)} $E_{123}^{12\ldots k, 12\ldots \ell}(x)=$
\[ \qquad=\left\{ \begin{array}{ll} 0, & \mbox{if $k \geq 3$ or $\ell \geq 3$,} \\ x-\frac{1}{2}-\frac{\sqrt{3}}{2}\tan\left( \frac{\sqrt{3}}{2}x + \frac{\pi}{6}\right)+\\ \sec\left( \frac{\sqrt{3}}{2}x + \frac{\pi}{6}\right)\left( \frac{\sqrt{3}}{2}\left( e^{x/2} + e^{-x/2} \right) - \sin\left( \frac{\sqrt{3}}{2}x + \frac{\pi}{3}\right) \right), & \mbox{if $k=2$ and $\ell = 2$,} \\ \frac{\sqrt{3}}{2}e^{x/2}\sec\left(\frac{\sqrt{3}}{2}x + \frac{\pi}{6}\right) - 1, & \mbox{if $k=1$ and $\ell = 1$,} \\ \frac{\sqrt{3}}{2}e^{x/2}\sec\left(\frac{\sqrt{3}}{2}x +
\frac{\pi}{6}\right) - \frac{1}{2} -
\frac{\sqrt{3}}{2}\tan\left(\frac{\sqrt{3}}{2}x +
\frac{\pi}{6}\right), & \mbox{else;}
\end{array}
\right.\]

{\rm(ii)} $E_{123}^{12\ldots k, \ell(\ell-1) \ldots 1}(x) =$
\[ \qquad=\left\{ \begin{array}{ll} 0, & \mbox{if $k \geq 3$,} \\ {\Phi}_{\ell}(x) = \frac{e^{x/2}}{(\ell-1)!}\sec\left(\frac{\sqrt{3}}{2}x + \frac{\pi}{6}\right)\int_{0}^{x}e^{-t/2}t^{\ell-1}\sin\left(\frac{\sqrt{3}}{2}t
+ \frac{\pi}{3}\right)\ dt, & \mbox{if $k = 1$,} \\ \int_{0}^{x}\sec\left(\frac{\sqrt{3}}{2}t + \frac{\pi}{6}\right)\left( \sin\left(\frac{\sqrt{3}}{2}t + \frac{\pi}{3}\right) - \frac{\sqrt{3}}{2}e^{-t/2}\right)\left( {\Phi}_{\ell}(t) + \frac{t^{\ell-1}}{(\ell-1)!} \right)\ dt, & \mbox{if $k = 2$;}
\end{array}
\right. \]

{\rm(iii)} $E_{123}^{k(k-1)\ldots 1, 12\ldots \ell}(x) =$
\[ \qquad=\left\{ \begin{array}{ll} 0, & \mbox{if $\ell \geq 3$,} \\ {\Phi}_k(x) = \frac{e^{x/2}}{(k-1)!}\sec\left(\frac{\sqrt{3}}{2}x + \frac{\pi}{6}\right)\int_{0}^{x}e^{-t/2}t^{k-1}\sin\left(\frac{\sqrt{3}}{2}t
+ \frac{\pi}{3}\right)\ dt, & \mbox{if $\ell = 1$,} \\ \int_{0}^{x}\sec\left(\frac{\sqrt{3}}{2}t + \frac{\pi}{6}\right)\left( \sin\left(\frac{\sqrt{3}}{2}t + \frac{\pi}{3}\right) - \frac{\sqrt{3}}{2}e^{-t/2}\right)\left( {\Phi}_k(t) + \frac{t^{k-1}}{(k-1)!} \right)\ dt, & \mbox{if $\ell = 2$;}
\end{array}
\right. \]

{\rm(iv)} $E_{123}^{k(k-1)\ldots 1, \ell(\ell-1)\ldots 1}(x)=\left\{ \begin{array}{ll} E_{123}^{\ell (\ell - 1) \ldots 1}(x), & \mbox{if $k = 1$,} \\ E_{123}^{k(k - 1) \ldots 1}(x), & \mbox{if $\ell = 1$,} \\ E_{123}^{k(k - 1) \ldots 1}(x) - E_{123}^{k(k - 1) \ldots 1,12}(x), & \mbox{if $\ell = 2$;} \\
\end{array}\right.$\\
For $k \geq 2$ and $\ell \geq 3$, $E_{123}^{k(k-1)\ldots 1,
\ell(\ell-1)\ldots 1}(x)$ satisfies
$$\frac{\partial}{\partial x}E_{123}^{k(k-1)\ldots 1,
\ell(\ell-1)\ldots 1}(x)= \left( E_{123}^{\ell(\ell - 1)\ldots
1}(x) + \frac{x^{\ell - 1}}{(\ell - 1)!}
\right)E_{123}^{k(k-1)\ldots 1,21}(x) + E_{123}^{(k-1) \ldots 1,
\ell(\ell-1)\ldots 1}(x),$$ where $E_{123}^{k(k-1)\ldots 1}(x)$ is
given in~\cite[Theorem 2]{KitMans}:
$$E_{123}^{k(k-1)\ldots 1}(x) = \frac{e^{x/2}\int_{0}^{x}e^{-t/2}t^{k-1}\sin(\frac{\sqrt{3}}{2}t + \frac{\pi}{6}))\ dt}{(k-1)!\cos(\frac{\sqrt{3}}{2}x + \frac{\pi}{6})},$$
and $E_{123}^{k(k-1)\ldots 1, 12}$ is given in this theorem above.
\end{prop}
\begin{proof}
\

{\bf Beginning with $k(k-1)\ldots 1$ and ending with $12\ldots
\ell$:} If $\ell \geq 3$ then the pattern $12\ldots \ell$ does not
avoid $123$, thus the statement is true. If $\ell = 1$, the
statement is true according to~\cite[Theorem 8]{Kit3} and the
observation that if $k=1$ then this formula gives the expression
$$\frac{\sqrt{3}}{2}e^{x/2}\sec\left(\frac{\sqrt{3}}{2}x + \frac{\pi}{6}\right) - 1,$$
which is true according to~\cite[Theorem 4.1]{ElizNoy} and the
assumption that the empty permutation does not begin or end with
the pattern $p=1$. So, we need only to consider the case $\ell =
2$. Recall the definitions of $E_{q}^{p}(x)$ in
Section~\ref{Preliminaries}.

Let $P_k(n)$ denote the number of $n$-permutations that avoid the
pattern $123$, begin with a decreasing subword of length $k$ and
end with the pattern $12$. Also, let $R_k(n)$ denote the number of
$n$-permutations that avoid the pattern $123$ and begin with a
decreasing subword of length $k$. Let $\pi = \pi_1 1 \pi_2$ be an
$(n+1)$-permutation that avoids the pattern $123$, begins with the
pattern $k(k-1)\ldots 1$ and ends with the pattern $12$. We
observe that $\pi_1$ avoids 123 and begins with $k(k-1)\ldots 1$;
$\pi_2$ ends with the pattern $12$ and $|\pi_2|>0$ since otherwise
$\pi$ cannot end with the pattern $12$; if $|\pi_2|>1$ then
$\pi_2$ must begin with the pattern $21$ since otherwise we have
an occurrence of the pattern $123$ beginning from the letter $1$.
If $|\pi_1|=i$ then the letters of $\pi_1$ can be chosen in
$\binom{n}{i}$ ways. So, there are at least
$$\sum_{i\geq0}\binom{n}{i}R_k(i)P_2(n-i) + nR_k(n-1)$$
$(n+1)$-permutations with the good properties, where the first
term corresponds to the case $|\pi_2|>1$ and the second term to
the case $|\pi_2|=1$. By this formula, we do not count the
permutations having $|\pi_1|=k-1$, although in this case $\pi$
begins with the pattern $k(k-1)\ldots 1$. So, we can choose the
letters of $\pi_1$ in $\binom{n}{k-1}$ ways, and according to
whether $|\pi|\geq 1$ or $|\pi|=1$, we have two terms:
$$\binom{n}{k-1}P_2(n-k+1) + k\delta_{n,k},$$
where $\delta_{n,k}$ is the Kronecker delta. Thus,
$$P_k(n+1) = \sum_{i\geq 0}\binom{n}{i}R_k(i)P_2(n-i) + nR_k(n-1) + \binom{n}{k-1}P_2(n-k+1) + k\delta_{n,k}.$$

After multiplying both sides of the last equality with $x^n/n!$
and summing over $n$, we have
\begin{equation}
\frac{d}{dx}E_{123}^{k(k-1)\ldots 1, 12}(x) = (E_{123}^{21, 12}(x) + x)\left(E_{123}^{k(k-1)\ldots 1}(x) + \frac{x^{k-1}}{(k-1)!}\right),
\label{eq15}
\end{equation}
with the initial condition $E_{123}^{k(k-1)\ldots 1, 12}(0)=0$. Since
$$E_{123}^{k(k-1)\ldots 1}(x)=E_{123}^{k(k-1)\ldots
1,1}(x)=\frac{e^{x/2}}{(k-1)!}\sec\left(\frac{\sqrt{3}}{2}x +
\frac{\pi}{6}\right)\int_{0}^{x}e^{-t/2}t^{k-1}\sin\left(\frac{\sqrt{3}}{2}t
+ \frac{\pi}{3}\right)\ dt,$$ to solve~(\ref{eq15}), we only need
to know $E_{123}^{21, 12}(x)$. To find it, we set $k=2$
into~(\ref{eq15}) and solve this equation. For an example how to
solve such an equation, we refer to~\cite[Theorem 6]{Kit3}. We get
$$E_{123}^{21, 12}(x)=-x+\sec\left( \frac{\sqrt{3}}{2}x + \frac{\pi}{6}\right)e^{-x/2}\int_{0}^{x}e^{t/2}\cos\left( \frac{\sqrt{3}}{2}t + \frac{\pi}{6}\right)\ dt.$$

Now, we put the formula for $E_{123}^{21,12}(x)$ into~(\ref{eq15})
and solve this differential equation to get the desired result.

{\bf Beginning with $12\ldots k$ and ending with
$\ell(\ell-1)\ldots 1$:} By the reverse and complement operations,
to avoid $123$, begin with the pattern $12\ldots k$ and end with
the pattern $\ell(\ell-1)\ldots 1$ is the same as to avoid $123$,
begin with the pattern $\ell(\ell-1)\ldots 1$ and end with the
pattern $12\ldots k$, so we can apply the results of the previous
case.

{\bf Beginning with $12\ldots k$ and ending with $12\ldots \ell$:}
The statement is obvious if $k \geq 3$ or $\ell \geq 3$. If $k=1$
and $\ell = 1$ then the statement is true according
to~\cite[Theorem 4.1]{ElizNoy} (but we need to subtract 1, since
by our assumption the empty permutation does not begin or end with
the pattern $p=1$). If $\ell = 1$ and $k=2$, the statement is true
according~\cite[Theorem 9]{Kit3}. If $k=1$ and $\ell = 2$, we
apply the reverse and complement operations, and use
again~\cite[Theorem 9]{Kit3}. So, we only need to consider the
case $k=2$ and $\ell = 2$. It is easy to see that
$$E_{123}^{12,12}(x) = E_{123}^{1,12}(x) - E_{123}^{21,12}(x),$$
and from the previous  cases
$$E_{123}^{1,12}(x) = \frac{\sqrt{3}}{2}e^{x/2}\sec\left( \frac{\sqrt{3}}{2}x + \frac{\pi}{6}\right) - \frac{1}{2} - \frac{\sqrt{3}}{2}\tan\left( \frac{\sqrt{3}}{2}x + \frac{\pi}{6}\right),$$
and $$E_{123}^{21,12}(x) = -x + \sec\left( \frac{\sqrt{3}}{2}x + \frac{\pi}{6}\right)\left(\sin\left( \frac{\sqrt{3}}{2}x + \frac{\pi}{3}\right)  - \frac{\sqrt{3}}{2}e^{-x/2}\right).$$

{\bf Beginning with $k(k-1)\ldots 1$ and ending with
$\ell(\ell-1)\ldots 1$:} If $\ell =1$, the statement is trivial.
If $k=1$, we get the statement by using the reverse and complement
operations. For the case $\ell = 2$, we observe that the number of
$n$-permutations that avoid the pattern $123$, begin with the
pattern $k(k-1)\ldots 1$ and end with the pattern $21$ is equal to
the number of $n$-permutation that avoid $123$ and begin with the
pattern $k(k-1)\ldots 1$ minus the number of $n$-permutations that
avoid the pattern $123$, begin with the pattern $k(k-1)\ldots 1$
and end with the pattern $12$. Suppose now that $k\geq 2$ and
$\ell \geq 3$ and an $(n+1)$-permutation $\pi$ avoids $123$,
begins with $k(k-1)\ldots 1$ and ends with $\ell(\ell-1)\ldots 1$.
It is easy to see that the letter $(n+1)$ can be either in the
first position, or in the position $i$, where $(k+1) \leq i \leq
(n-\ell)$, or in the position $(n- \ell +1)$. In the first of
these cases, obviously we have $N_{123}^{(k-1)\ldots 1, \ell (\ell
- 1) \ldots 1}(n)$ permutations. In the second case, we choose the
letters of $\pi$ to the left of $(n+1)$ in $\binom{n}{i-1}$ ways.
These letters must form a permutation that avoids $123$, begins
with the pattern $k(k-1)\ldots 1$, and ends with the pattern $21$
(if the last condition does not hold, the letter $(n+1)$ and two
letters to the left of it form a $123$-pattern. In the same time,
the letters to the right of $(n+1)$ form a permutation that avoids
$123$ and ends with the pattern $\ell(\ell-1)\ldots 1$. In the
third case, we can choose the letters to the right of $(n+1)$ in
$\binom{n}{\ell-1}$ ways, rearrange them into the decreasing
order, and form from the letters to the left of $(n+1)$ a
permutation that avoids $123$, begins with the pattern
$k(k-1)\ldots 1$ and ends with the pattern $21$ (by the same
reasons as above) in $N_{123}^{k(k-1)\ldots 1,21}(n-\ell + 1)$
ways. Thus,
$$\begin{array}{l}
N_{123}^{k(k-1)\ldots 1,\ell(\ell-1)\ldots
1}(n+1)=N_{123}^{(k-1)\ldots 1, \ell (\ell - 1) \ldots 1}(n)\\
\qquad\qquad\qquad+\displaystyle\sum_{i=0}^{n}\binom{n}{i}N_{123}^{k(k-1)\ldots
1,21}(i)N_{123}^{\ell (\ell -1)\ldots 1}(n-i) +
\binom{n}{\ell-1}N_{123}^{k(k-1)\ldots 1,21}(n-\ell +
1),\end{array}$$ where we observed, that to avoid $123$ and end
with $\ell (\ell -1)\ldots 1$ is the same as to avoid $123$ and
begin with $\ell (\ell -1)\ldots 1$ using the reverse and
complement. Now, we multiply both sides of the equality by
$x^n/n!$ and sum over $n$ to get the desirable result.
\end{proof}

\section{Avoiding a pattern x-yz, beginning and ending with certain patterns simultaneously}

\begin{prop}\label{prop16} We have

{\rm(i)}
$E_{1\mn32}^{12\ldots k, 1}(x)= E_{1\mn32}^{12\ldots k}(x) = \left\{ \begin{array}{ll} e^{e^x}\int_{0}^{x}e^{-e^t}\displaystyle\sum_{n\geq k-1}\frac{t^n}{n!}\ dt, & \mbox{if $k\geq 2$}\\[2mm] e^{e^x-1}, & \mbox{if $k=1$}
\end{array}\right.$.\\ For $\ell \geq 2$, $E_{1\mn32}^{12\ldots k, 12\ldots
\ell}(x)$ satisfies
$$\frac{\partial}{\partial x}E_{1\mn32}^{12\ldots k, 12\ldots \ell}(x)= \left(e^x - \displaystyle\sum_{i=0}^{\ell-2}\frac{x^i}{i!}\right)E_{1\mn32}^{12\ldots k}(x) + e^xx^{max(\ell,k)-1}.$$

{\rm(ii)} $E_{1\mn32}^{12\ldots k, \ell(\ell-1) \ldots 1}(x)$
satisfies
$$\frac{{\partial}^{\ell - 1}}{\partial x^{\ell - 1}}E_{1\mn32}^{12\ldots k, \ell(\ell-1) \ldots 1}(x)=\left\{ \begin{array}{ll} e^{e^x}\int_{0}^{x}e^{-e^t}\displaystyle\sum_{n\geq k-1}\frac{t^n}{n!}\ dt, & \mbox{if $k\geq 2$,}\\[2mm] e^{e^x-1}, & \mbox{if $k=1$.}
\end{array}
\right.$$

{\rm(iii)}
$E_{1\mn32}^{k(k-1)\ldots 1, 1}(x)= E_{1\mn32}^{k(k-1)\ldots 1}(x) = \left\{ \begin{array}{ll} (e^{e^x}/(k-1)!)\int_{0}^{x}t^{k-1}e^{-e^t+t}\ dt, & \mbox{if $k\geq 2$}\\[2mm] e^{e^x-1}, & \mbox{if $k=1$}
\end{array}\right.$.\\ For $\ell \geq 2$, $E_{1\mn32}^{k(k-1)\ldots 1,12\ldots \ell}(x)$ satisfies
$$\frac{\partial}{\partial x}E_{1\mn32}^{k(k-1)\ldots 1, 12\ldots \ell}(x)=\left(e^x - \displaystyle\sum_{i=0}^{\ell-2}\frac{x^i}{i!}\right)E_{1\mn32}^{k(k-1)\ldots 1}(x) + \left( e^x - \displaystyle\sum_{i=0}^{\ell-2}\frac{x^i}{i!}\right)\frac{x^{k-1}}{(k-1)!}.$$

{\rm(iv)} $E_{1\mn32}^{k(k-1)\ldots 1, \ell(\ell-1) \ldots 1}(x)$
satisfies
$$
\frac{{\partial}^{\ell - 1}}{\partial x^{\ell -
1}}\left(E_{1\mn32}^{k(k-1)\ldots 1, \ell(\ell-1) \ldots 1}(x) -
\frac{x^{max(k,\ell)}-x^{k+\ell-1}}{1-x}\right)=\left\{ \begin{array}{ll} \frac{e^{e^x}}{(k-1)!}\int_{0}^{x}t^{k-1}e^{-e^t+t}\ dt, & \mbox{if $k\geq 2$,}\\[2mm] e^{e^x-1}, & \mbox{if $k=1$.}
\end{array}\right.$$
\end{prop}
\begin{proof}
\

{\bf Beginning with $12\ldots k$ and ending with
$\ell(\ell-1)\ldots 1$:} If $\ell = 1$ then the result follows
from \cite[Proposition 5]{KitMans}, since to avoid $1\mn32$ and
begin with $12\ldots k$ is the same as to avoid $3\mn12$ and begin
with $k(k-1)\ldots 1$. Suppose now that $\ell \geq 2$ and a
permutation $\pi$ avoids the pattern $1\mn32$, begins with the
pattern $12\ldots k$ and ends with the pattern $\ell(\ell-1)\ldots
1$. Since $\ell \geq 2$, we have that the letter 1 must be in the
rightmost position since otherwise, this letter and two rightmost
letters of $\pi$ form the pattern $1\mn32$, which is forbidden.
Thus,
$$N^{12\ldots k,\ell(\ell-1)\ldots 1}_{1\mn32}(n)=N^{12\ldots k,(\ell-1)(\ell-2)\ldots 1}_{1\mn32}(n-1)=\cdots =N^{12\ldots k, 1}_{1\mn32}(n-\ell+1).$$
Multiplying both sides of the equality $N^{12\ldots
k,\ell(\ell-1)\ldots 1}_{1\mn32}(n)=N^{12\ldots k,
1}_{1\mn32}(n-\ell+1)$ by $x^{n-\ell+1}/(n-\ell+1)!$ and summing
over $n$, we get
$$\frac{{\partial}^{\ell - 1}}{\partial x^{\ell - 1}}E_{1\mn32}^{12\ldots k, \ell(\ell-1) \ldots 1}(x)=E_{1\mn32}^{12\ldots k}(x),$$
where $E_{1\mn32}^{12\ldots k}(x)$ is given in~\cite[Proposition
5]{KitMans}, since to avoid $1\mn32$ and begin with $12\ldots k$
is the same as to avoid $3\mn12$ and begin with $k(k-1)\ldots 1$.

{\bf Beginning with $k(k-1)\ldots 1$ and ending with
$\ell(\ell-1)\ldots 1$:} We use the same arguments as those given
under consideration of the previous case, but instead
of~\cite[Proposition 5]{KitMans} we use~\cite[Proposition
4]{KitMans}. However, we observe, that when we use the argument
$$N^{k(k-1)\ldots 1,\ell(\ell-1)\ldots 1}_{1\mn32}(n)=N^{k(k-1)\ldots 1,(\ell-1)(\ell-2)\ldots 1}_{1\mn32}(n-1)=\cdots =N^{k(k-1)\ldots 1, 1}_{1\mn32}(n-\ell+1)$$
for $k,\ell \geq 2$, we do not count the decreasing permutations
of length $max(k,\ell)$, $max(k,\ell)+1, \ldots , k+\ell-2$, since
in this case, the patterns $k(k-1)\ldots 1$ and
$\ell(\ell-1)\ldots 1$ overlap in more than one letter, which
causes the observation. So, we need to consider additionally the
term
$$x^{max(k,\ell)}+x^{max(k,\ell)+1}+\cdots+x^{k+\ell-2}=\frac{x^{max(k,\ell)}-x^{k+\ell-1}}{1-x},$$
which vanishes if $k=1$ or $\ell=1$.

{\bf Beginning with $12\ldots k$ and ending with $12\ldots \ell$:}
The only interesting case here is the case $k\geq 2$ and $\ell
\geq 2$. Using the reverse and complement, instead of considering
avoiding $1\mn32$, beginning with $12\ldots k$ and ending with
$12\ldots \ell$, we consider avoiding $1\mn32$, beginning with
$12\ldots \ell$ and ending with $12\ldots k$. Suppose an
$n$-permutation $\pi$ satisfies all the conditions. We observe,
that the letter $n$ can be in the position $i$, where $\ell \leq i
\leq n-k$. Also, $n$ can be in the rightmost position if $n \geq
max(\ell, k)$. In any case, the letters of $\pi$ to the left of
$n$ must be in the increasing order, since otherwise we have an
occurrence of the pattern $21\mn3$. This means that in the second
case we have the only one permutation. In the first case, the
letters of $\pi$ to the right of $n$ must avoid $21\mn3$ and end
with the pattern $12\ldots k$. The number of such permutations,
using the reverse and complement, is given by
$N_{1\mn32}^{12\ldots k}(n-i)$. Thus, for $n\geq max(\ell, k)$,
$$N_{21\mn3}^{12\ldots \ell, 12\ldots k}(n) = \displaystyle\sum_{i=\ell}^{n-k}\binom{n-1}{i-1}N_{1\mn32}^{12\ldots k}(n-i) + 1.$$
This gives
$$N_{21\mn3}^{12\ldots \ell, 12\ldots k}(n) = \displaystyle\sum_{i=1}^{n}\binom{n-1}{i-1}N_{1\mn32}^{12\ldots k}(n-i) - \displaystyle\sum_{i=1}^{\ell-1}\binom{n-1}{i-1}N_{1\mn32}^{12\ldots k}(n-i)+ 1,$$
which leads to the desirable result after multiplying both sides
of the last equality by $x^n/n!$ and summing over $n$.

{\bf Beginning with $k(k-1)\ldots 1$ and ending with $12\ldots
\ell$:} The only interesting case here is the case $k\geq 2$ and
$\ell \geq 2$. Using the reverse and complement, instead of
considering avoiding $1\mn32$, beginning with $k(k-1)\ldots 1$ and
ending with $12\ldots \ell$, we consider avoiding $1\mn32$,
beginning with $12\ldots \ell$ and ending with $k(k-1)\ldots 1$.
Suppose an $n$-permutation $\pi$ satisfies all the conditions. We
observe, that the letter $n$ can only be in the position $i$,
where $\ell \leq i \leq n-k$, or in position $(n-k+1)$ (in the
case $n\geq k+\ell-1$). In the first case, it is easy to see that
the letters of $\pi$ to the left of $n$ must be in the increasing
order, and the letters of $\pi$ to the right of $n$ must avoid
$21\mn3$ and end with the pattern $k(k-1)\ldots 1$. Using the
reverse and complement, the total number of permutations counted
in the first case is
$\displaystyle\sum_{i=\ell}^{n-k}\binom{n-1}{i-1}N_{1\mn32}^{k(k-1)\ldots
1}(n-i)$. In the second case, the letters to the left of $n$ are
in the increasing order, whereas the letters to the right of $n$
are in decreasing order. The number of such permutations is
$\binom{n-1}{k-1}$, which is the number of ways to choose the last
$k-1$ letters. Thus,
$$N_{21\mn3}^{12\ldots \ell, k(k-1)\ldots 1}(n) = \displaystyle\sum_{i=\ell}^{n-k}\binom{n-1}{i-1}N_{1\mn32}^{k(k-1)\ldots 1}(n-i) + \binom{n-1}{k-1}.$$
Multiplying both parts of the equality by $x^{n-1}/(n-1)!$ and summing over $n$, we get
$$\begin{array}{l}
\dfrac{\partial}{\partial x}E_{21\mn3}^{12\ldots \ell,
k(k-1)\ldots 1}(x)=\displaystyle\sum_{n\geq k+\ell}\binom{n-1}{k-1}\frac{x^{n-1}}{(n-1)!}\\
\qquad\qquad\qquad+\displaystyle\sum_{n\geq 0}\left(
\displaystyle\sum_{i=1}^{n-1}\binom{n-1}{i-1}N_{1\mn32}^{k(k-1)\ldots
1}(n-i) -
\displaystyle\sum_{i=1}^{\ell-1}\binom{n-1}{i-1}N_{1\mn32}^{k(k-1)\ldots
1}(n-i)\right)\frac{x^{n-1}}{(n-1)!},\end{array}$$ which leads to
the desirable result.
\end{proof}

\begin{prop}\label{prop17} We have

{\rm(i)} $G_{2\mn13}^{12\ldots k, 12\ldots \ell}(x) =  x^{k+\ell
-1}C^{k+1}(x)+\frac{x^m-x^{k+\ell -1}}{1-x}$.

{\rm(ii)} $G_{2\mn 13}^{k(k-1)\ldots 1, 12\ldots \ell}(x) =
x^{k+\ell -1}C^2(x)$.

{\rm(iii)} $G_{2\mn13}^{k(k-1)\ldots 1, \ell(\ell-1)\ldots 1}(x) =
x^{k+\ell -1}C^{\ell+1}(x)+\frac{x^m-x^{k+\ell -1}}{1-x}$, where
$m=max(k,\ell)$.

{\rm(iv)} the generating function $G_{2\mn13}(x,y,z)=\sum_{k,\ell
\geq 0} G_{2\mn13}^{12\ldots k, \ell(\ell-1) \ldots
1}(x)y^kz^{\ell}$ for the sequence\\ $\{G_{2\mn13}^{12\ldots k,
\ell(\ell-1) \ldots 1}(x)\}_{k,\ell\geq0}$ {\rm(}where $k$ and
$\ell$ go through all natural numbers{\rm)} is
$$\frac{1}{1-x(y+z)}\left( x(y+z+yz) + \frac{C(x)-1}{(1-xyC(x))(1-xzC(x))}\right).$$
\end{prop}
\begin{proof}
\ By \cite{Claes}, to avoid the pattern $2\mn13$ is the same as to
avoid the pattern $2\mn1\mn3$. Thus we can apply the results of
Proposition~\ref{prop11}.
\end{proof}

\begin{prop}\label{prop18} We have

{\rm(i)} $E_{1\mn23}^{12\ldots k, 12\ldots \ell}(x) = \left \{
\begin{array}{ll} 0, & \mbox{if $k\geq 3$ or $\ell \geq 3$,} \\
E_{1\mn23}^{12\ldots k}(x), & \mbox{if $\ell = 1$,} \\
E_{12\mn3}^{12\ldots \ell}(x), & \mbox{if $k = 1$,} \\
\int_{0}^{x}tE_{12\mn3}^{12}(t)\ dt + \frac{x^2}{2!}, &
\mbox{if $k=2$ and $\ell = 2$,} \\
\end{array}\right.$\\
where $E_{12\mn3}^{12\ldots k}(x)$ and $E_{1\mn23}^{12\ldots
k}(x)$ are given by~\cite[Proposition 10]{KitMans}
and~\cite[Proposition 6]{KitMans} respectively:

$E^{12\dots k}_{12\mn3}(x)=\left\{\begin{array}{ll}
0, & \mbox{if $k\geq3$,} \\[2mm]
x^2\displaystyle\sum_{j=0}^k(1-jx)^{-1}\displaystyle\sum_{d\geq0}\frac{x^d}{(1-x)(1-2x)\dots(1-dx)}, & \mbox{if $k=2$,} \\[4mm]
\displaystyle\sum_{d\geq 0}\frac{x^d}{(1-x)(1-2x)\dots(1-dx)}, & \mbox{if $k=1$;}
\end{array}\right.$

$E_{1\mn23}^{12\ldots k}(x) = E_{3-21}^{k(k-1)\ldots 1}(x) = \left\{ \begin{array}{ll} 0, & \mbox{if $k\geq 3$,} \\[2mm] e^{e^x}\int_{0}^{x}e^{-e^t}(e^t - 1)\ dt, & \mbox{if $k=2$,}
\\[2mm] e^{e^x-1}, & \mbox{if $k=1$.}
\end{array}\right.$

{\rm(ii)} $N_{1\mn23}^{12\ldots k, \ell(\ell-1) \ldots 1}(n)=$
\[\quad\left
\{ \begin{array}{ll} 0, & \mbox{if $k \geq 3$,} \\ 0, & \mbox{if $k=2$ and $n\leq \ell$,} \\ 1 + N_{1\mn23}^{12, (\ell-1)(\ell-2)\ldots 1}(n-1) + \displaystyle\sum_{j=\ell+1}^{n-2}\binom{n-1}{j-1}N_{1\mn23}^{12}(n-j), & \mbox{if $k=2$ and $n\geq \ell+1$,} \\ N_{12\mn3}^{\ell(\ell-1)\ldots 1}(n), & \mbox{if $k=1$,} \\
\end{array}\right.\] where the numbers $N_{12\mn3}^{\ell (\ell-1) \ldots 1}(n)$
are given in~\cite[Proposition 9]{KitMans}, and the numbers
$N_{1\mn23}^{12}(n)$ are given by expending the exponential
generating functions in~\cite[Proposition~6]{KitMans}.

{\rm(iii)} $E_{1\mn23}^{k(k-1)\ldots 1, 12\ldots \ell}(x) = \left \{ \begin{array}{ll} 0, & \mbox{if $\ell \geq 3$} \\[2mm] \frac{1}{(k-1)!}\int_{0}^{x}\int_{0}^{t}tm^{k-1}e^{e^t-e^m+m}\ dmdt+\frac{kx^{k+1}}{(k+1)!}, & \mbox{if $\ell=2$} \\[2mm] (e^{e^x}/(k-1)!)\int_{0}^{x}t^{k-1}e^{-e^t+t}\ dt, & \mbox{if $\ell=1$} \\
\end{array}\right.$,\\ where
$E_{1\mn23}^{k(k-1)\ldots 1,1}(n) = E_{1\mn23}^{k(k-1)\ldots
1}(n)$ is given by~\cite[Proposition 4]{KitMans}, and
$N_{1\mn23}^{1, \ell(\ell-1)\ldots 1}(n) =
N_{12\mn3}^{\ell(\ell-1)\ldots 1}(n)$ is given
by~\cite[Proposition 9]{KitMans};

{\rm(iv)} For $k\geq 2$ and $\ell \geq 2$,
$E_{1\mn23}^{k(k-1)\ldots 1, \ell(\ell-1)\ldots 1}(x)$ satisfies
$$\frac{\partial}{\partial x}E_{1\mn23}^{k(k-1)\ldots 1,
\ell(\ell-1)\ldots 1}(x) = E_{1\mn23}^{k(k-1)\ldots 1,
(\ell-1)\ldots 1}(x)+\left( e^x -
\displaystyle\sum_{i=0}^{\ell-1}\frac{x^i}{i!}\right)\left(
E_{1\mn23}^{k(k-1)\ldots 1}(x) + \frac{x^{k}}{(k-1)!}\right).$$
\end{prop}
\begin{proof}
\

{\bf Beginning with $k(k-1)\ldots 1$ and ending with $12\ldots
\ell$:} If $\ell \geq 3$ then $E_{1\mn23}^{k(k-1)\ldots 1,
12\ldots \ell}(x)=0$, since in this case the pattern $12\ldots
\ell$ does not avoid $1\mn23$. If $\ell = 1$ then we
use~\cite[Proposition 4]{KitMans}, since in this case the only
restrictions to the permutations are avoiding $1\mn23$ and
beginning with the pattern $k(k-1)\ldots 1$. Suppose now that
$\ell = 2$ and an $(n+1)$-permutation $\pi$ avoids $1\mn23$,
begins with $k(k-1)\ldots 1$ and ends with the pattern $12$. The
letter $1$ must be in next to the rightmost position, since
otherwise this letter and two rightmost letters form the pattern
$1\mn23$. We can choose the rightmost letter of $\pi$ in $n$ ways,
and the letters to the left of $1$ must form a $1\mn23$-avoiding
permutation that begins with $k(k-1)\ldots 1$. Besides, if $n=k$,
and the $k-1$ letters to the right of $1$ are in the decreasing
order, we get $n$ extra permutations that satisfy our
restrictions. Thus,
$$N_{1\mn23}^{k(k-1)\ldots 1, 12}(n+1) = nN_{1\mn23}^{k(k-1)\ldots 1}(n) + n{\delta}_{n,k},$$
where $\delta_{n,k}$ is the Kronecker delta. Multiplying both
sides of the equality by $x^n/n!$ and summing over $n$ we get
$$E_{1\mn23}^{k(k-1)\ldots 1, 12}(x)=\int_{0}^{x}tE_{1\mn23}^{k(k-1)\ldots 1}(t)\ dt+\frac{kx^{k+1}}{(k+1)!}.$$
Using the formula for $E_{1\mn23}^{k(k-1)\ldots 1}(t)$
in~\cite[Proposition 4]{KitMans}, we get the desirable result.

{\bf Beginning with $12\ldots k$ and ending with $12\ldots \ell$:}
The first three cases are easy to prove in the same manner as we
do in the proves of previous propositions. The only interesting
case is when $k=2$ and $\ell =2$. Using the reverse and complement
operations, instead of considering avoiding $1\mn23$, beginning
with $12$ and ending with $12$, we consider avoiding $12\mn3$,
beginning with $12$ and ending with $12$, which we find to be more
easy. Suppose an $(n+1)$-permutation $\pi$ satisfies all the
restrictions. It is easy to see that $|\pi| \neq 1$ and $|\pi|
\neq 3$, as well as if $|\pi|=2$ (that is $n=1$) then $\pi$ must
be $12$. Suppose $|\pi| \geq 4$. Since $\pi$ begins with the
pattern $12$, it is impossible for the letter $(n+1)$ to be
somewhere to the right of the second letter of $\pi$ or to be the
leftmost letter. Thus, $(n+1)$ must be in the second position. We
can choose the leftmost letter of $\pi$ in $n$ ways, since any
choice of this letter will not lead to an occurrence of the
pattern $12-3$ beginning with two leftmost letters. If $\pi =
a(n+1){\pi}^{\prime}$ then ${\pi}^{\prime}$ must avoid $12\mn3$
and end with the pattern $12$. The number of such permutations,
using the reverse and complement, is given by
$N_{1\mn23}^{12}(n-1)$. Thus,
$$N_{12\mn3}^{12,12}(n+1)=nN_{1\mn23}^{12}(n-1).$$
Multiplying both sides of the equality by $x^n/n!$ and summing over all $n$, we get
$$(E_{12\mn3}^{12,12}(x))^{\prime}=xE_{1\mn23}^{12}(x)+ x,$$
where the term $x$ corresponds to the permutation $12$. We have
the desirable result by integrating both sides of the last
equality.

{\bf Beginning with $12\ldots k$ and ending with
$\ell(\ell-1)\ldots 1$:} All the cases but $k=2$ and $n \geq \ell
+ 1$ are easy to prove. Let us consider this case. Using the
reverse and complement operations, instead of considering avoiding
$1\mn23$, beginning with $12$ and ending with $\ell(\ell -1)
\ldots 1$, we consider avoiding $12\mn3$, beginning with
$\ell(\ell -1) \ldots 1$ and ending with $12$, which we find to be
more easy. Let an $n$-permutation $\pi$ satisfies all the
conditions. We observe, that the letter $n$ is either in the first
position, or in position $j$, where $k+1 \leq j \leq n-2$, or in
the last position. Obviously, in the first of these cases the
number of ``good'' permutations is given by $N_{12\mn3}^{(\ell
-1)(\ell - 2) \ldots 1, 12}(n-1)$, which is equivalent to
$N_{1\mn23}^{12, (\ell -1)(\ell - 2) \ldots 1}(n-1)$ by using the
reverse and complement. In the second case, we choose the letters
to the left of $n$ in $\binom{n-1}{j-1}$ ways, rearrange them to
the decreasing order (we do it since otherwise we have an
occurrence of the pattern $12\mn3$ having the letter $n$). After
that, the letters to the right of $n$ must form a permutation that
avoid $12\mn3$ and end with the pattern $12$. Using the reverse
and complement, there are $N_{1\mn23}^{12}(n-j)$ such
permutations. So, totally, in the second case there are
$\sum_{j=\ell+1}^{n-2}\binom{n-1}{j-1}N_{1\mn23}^{12}(n-j)$
permutations. Finally, if $n$ is at the last position, we have the
only one such permutation, since the other letters must be in the
decreasing order.

{\bf Beginning with $k(k-1)\ldots 1$ and ending with
$\ell(\ell-1)\ldots 1$:} The only interesting case here is the
case $k\geq 2$ and $\ell \geq 2$. Using the reverse and complement
operations, instead of considering avoiding $1\mn23$, beginning
with $k(k-1)\ldots 1$ and ending with $\ell(\ell -1) \ldots 1$, we
consider avoiding $12\mn3$, beginning with $\ell(\ell -1) \ldots
1$ and ending with $k(k-1)\ldots 1$, which we find to be more
easy. Let an $n$-permutation $\pi$ satisfies all the conditions.
We observe, that the letter $n$ is either in the first position,
or in position $j$, where $\ell+1 \leq j \leq n-k$, or in the last
position $n-k+1$. We proceed as in the previous case to get the
following
$$N_{12\mn3}^{\ell(\ell -1) \ldots 1, k(k-1)\ldots 1}=N_{12\mn3}^{(\ell -1) \ldots 1, k(k-1)\ldots 1} + \displaystyle\sum_{i=\ell+1}^{n-k}\binom{n-1}{i-1}N_{1\mn23}^{k(k-1)\ldots 1}(n-i) + \binom{n-1}{k-1},$$
where three terms in the right-hand side correspond to the three
cases described above. We now multiply both sides of the equality
by $x^n/n!$, sum over $n$ and observe the following detail. We
cannot write instead of $i=\ell+1$ (in the sum above) $i=1$ as we
did in most of the cases above, since, for instance, the case
$i=1$ do not necessarily make the term of summation equal 0 as it
was before. Thus, instead of the factor $e^x$, we have the factor
$\left( e^x- \displaystyle\sum_{i=0}^{\ell-1}\frac{x^i}{i!}
\right)$
\end{proof}
\section{Avoiding a pattern xy-z, beginning and ending with certain patterns simultaneously}

\begin{prop} We have

{\rm(i)} $G_{13\mn2}^{12\ldots k, 12\ldots \ell}(x) = x^{k+\ell
-1}C^{\ell+1}(x)+\frac{x^m-x^{k+\ell -1}}{1-x}$.

{\rm(ii)} $G_{13\mn2}^{12\ldots k, \ell(\ell-1) \ldots 1}(x) =
x^{k+\ell-1}C^2(x)$.

{\rm(iii)} $G_{13\mn2}^{k(k-1)\ldots 1, \ell(\ell-1)\ldots 1}(x) =
x^{k+\ell-1}C^{k+1}(x)+\frac{x^m-x^{k+\ell -1}}{1-x}$, where
$m=max(k,\ell)$.

{\rm(iv)} the generating function $G_{13-2}(x,y,z)=\sum_{k,\ell
\geq 0}G_{13\mn2}^{k(k-1)\ldots 1, 12\ldots \ell}(z)y^kz^{\ell}$
for the sequence\\ $\{G_{13\mn2}^{k(k-1)\ldots 1, 12\ldots
\ell}(x)\}_{k,\ell\geq0}$ {\rm(}where $k$ and $\ell$ go through
all natural numbers{\rm)} is
$$\frac{1}{1-x(y+z)}\left( x(y+z+yz) + \frac{C(x)-1}{(1-xyC(x))(1-xzC(x))}\right).$$
\end{prop}
\begin{proof}
We apply the reverse and complement operations and then use the
results of Proposition~\ref{prop17}. For example, to avoid
$2\mn13$, begin with $12\ldots k$ and end with $12\ldots \ell$ is
the same as to avoid $13\mn2$, begin with $12\ldots \ell$ and end
with $12\ldots k$.
\end{proof}

\begin{prop}\label{prop20} We have

{\rm(i)} $E_{21\mn3}^{12\ldots k, 1}(x)= E_{21\mn3}^{12\ldots
k}(x)$ is given by~\cite[Proposition~14]{KitMans}. For $\ell \geq
2$, $E_{21\mn3}^{12\ldots k, 12\ldots \ell}(x)$ satisfies
$$\frac{\partial}{\partial x}E_{21\mn3}^{12\ldots k, 12\ldots \ell}(x)= \left(e^x - \displaystyle\sum_{i=0}^{k-2}\frac{x^i}{i!}\right)E_{1\mn32}^{12\ldots \ell}(x) + e^xx^{max(\ell,k)-1}.$$
where $E_{1\mn32}^{12\ldots \ell}(x)=E_{3\mn12}^{\ell (\ell
-1)\ldots 1}(x)$ is given by~\cite[Proposition 5]{KitMans}.

{\rm(ii)} For $\ell \geq 2$, $E_{21\mn3}^{12\ldots k,
\ell(\ell-1)\ldots 1}(x)$ satisfies
$$\frac{\partial}{\partial x}E_{21\mn3}^{12\ldots k, \ell(\ell-1)\ldots 1}(x)=\left(e^x - \displaystyle\sum_{i=0}^{k-2}\frac{x^i}{i!}\right)E_{1\mn32}^{\ell(\ell-1)\ldots 1}(x) + \left( e^x - \displaystyle\sum_{i=0}^{k-2}\frac{x^i}{i!}\right)\frac{x^{\ell-1}}{(\ell-1)!},$$
where $E_{1\mn32}^{\ell(\ell-1)\ldots 1}(x)$ is given
by~\cite[Proposition 4]{KitMans}.

{\rm(iii)} $E_{21\mn3}^{k(k-1)\ldots 1, 12\ldots \ell}(x)$
satisfies
$$\frac{{\partial}^{k - 1}}{\partial x^{k - 1}}E_{21\mn3}^{k(k-1)\ldots 1, 12\ldots \ell}(x)=\left\{ \begin{array}{ll} e^{e^x}\int_{0}^{x}e^{-e^t}\displaystyle\sum_{n\geq \ell-1}\frac{t^n}{n!}\ dt, & \mbox{if $\ell \geq 2$,}\\[2mm] e^{e^x-1}, & \mbox{if $\ell=1$.}
\end{array}\right.$$

{\rm(iv)} $E_{21\mn3}^{k(k-1)\ldots 1, \ell(\ell-1) \ldots 1}(x)$
satisfies
$$\frac{{\partial}^{k - 1}}{\partial
x^{k-1}}\left(E_{2\mn13}^{k(k-1)\ldots 1, \ell(\ell-1) \ldots
1}(x)- \frac{x^{max(k,\ell)}-x^{k+\ell-1}}{1-x}\right)=\left\{ \begin{array}{ll} \frac{e^{e^x}}{(\ell-1)!}\int_{0}^{x}t^{\ell-1}e^{-e^t+t}\ dt, & \mbox{if $\ell \geq 2$,}\\[2mm] e^{e^x-1}, & \mbox{if $\ell=1$.}
\end{array}\right.$$
\end{prop}

\begin{proof}
We apply the reverse and complement operations and then use the
results of Proposition~\ref{prop16}. For example, to avoid
$1\mn32$, begin with $12\ldots k$ and end with $12\ldots \ell$ is
the same as to avoid $21\mn3$, begin with $12\ldots \ell$ and end
with $12\ldots k$.
\end{proof}

\begin{prop} We have

{\rm(i)} $E_{12\mn3}^{12\ldots k, 12\ldots \ell}(x) = \left \{
\begin{array}{ll} 0, & \mbox{if $k\geq 3$ or $\ell \geq 3$,} \\
E_{12\mn3}^{12\ldots k}(x), & \mbox{if $\ell = 1$,} \\
E_{1\mn23}^{12\ldots \ell}(x), & \mbox{if $k = 1$,} \\
\int_{0}^{x}tE_{1\mn23}^{12}(t)\ dt + \frac{x^2}{2!}, &
\mbox{if $k=2$ and $\ell = 2$,} \\
\end{array}
\right.$\\ where $E_{12\mn3}^{12\ldots k}(x)$ and
$E_{1\mn23}^{12\ldots k}(x)$ are given in
Proposition~\ref{prop18}.

{\rm(ii)} $E_{1\mn23}^{12\ldots k, \ell(\ell-1) \ldots 1}(x) =
\left
\{ \begin{array}{ll} 0, & \mbox{if $k \geq 3$,} \\[2mm] \frac{1}{(\ell-1)!}\int_{0}^{x}\int_{0}^{t}tm^{\ell-1}e^{e^t-e^m+m}\ dmdt+\frac{\ell x^{\ell+1}}{(\ell+1)!}, & \mbox{if $k=2$,} \\[2mm] (e^{e^x}/(\ell-1)!)\int_{0}^{x}t^{\ell-1}e^{-e^t+t}\ dt, & \mbox{if $k=1$;} \\
\end{array}\right. $

{\rm(iv)} $N_{12\mn3}^{k(k-1)\ldots 1, 12\ldots \ell}(n) =$
\[\left
\{ \begin{array}{ll} 0, & \mbox{if $\ell \geq 3$,} \\ 0, & \mbox{if $\ell=2$ and $n\leq k$,} \\ 1 + N_{12\mn3}^{(k-1)(k-2)\ldots 1, 12}(n-1) + \displaystyle\sum_{j=k+1}^{n-2}\binom{n-1}{j-1}N_{3\mn21}^{21}(n-j), & \mbox{if $\ell=2$ and $n\geq k+1$,} \\ N_{12\mn3}^{k(k-1)\ldots 1}(n), & \mbox{if $\ell=1$,} \\
\end{array}\right.\] where the numbers $N_{12\mn3}^{k(k-1)\ldots 1}(n)$ are
given in~\cite[Proposition 9]{KitMans}, and the numbers
$N_{3\mn21}^{21}(n)$ are given by expending the exponential
generating functions in~\cite[Proposition 6]{KitMans}.

{\rm(iv)} $N_{12\mn3}^{k(k-1)\ldots 1,1}(n) =
N_{12\mn3}^{k(k-1)\ldots 1}(n)$ is given by~\cite[Proposition
9]{KitMans}, and $N_{12\mn3}^{1, \ell(\ell-1)\ldots 1}(n) =
N_{1\mn23}^{\ell(\ell-1)\ldots 1}(n)$ is given
by~\cite[Proposition 4]{KitMans}. For $k\geq 2$ and $\ell \geq 2$,
$E_{12\mn3}^{k(k-1)\ldots 1, \ell(\ell-1)\ldots 1}(x)$ satisfies
$$\frac{\partial}{\partial x}E_{12\mn3}^{k(k-1)\ldots 1,
\ell(\ell-1)\ldots 1}(x) = E_{12\mn3}^{(k-1)\ldots 1,
\ell(\ell-1)\ldots 1}(x)+\left( e^x -
\displaystyle\sum_{i=0}^{k-1}\frac{x^i}{i!}\right)\left(
E_{1-23}^{\ell (\ell-1)\ldots 1}(x) +
\frac{x^{\ell}}{(\ell-1)!}\right).$$
\end{prop}

\begin{proof}
We apply the reverse and complement operations and then use the
results of Proposition~\ref{prop18}. For example, to avoid
$1\mn23$, begin with $12\ldots k$ and end with $12\ldots \ell$ is
the same as to avoid $12\mn3$, begin with $12\ldots \ell$ and end
with $12\ldots k$.
\end{proof}
\section{Further results} \label{section9}
In this section, we propose two directions of generalization of
the results from the previous sections. The first one is a
consideration of avoiding more than one pattern, beginning with
some pattern and ending with another pattern. For example, suppose
that $v=12\mn3$, $w=21\mn3$, $p=12\ldots k$, $q=12\ldots\ell$, and
$E_{v,w}^{p,q}(x)$ denotes the exponential generating function for
the number of permutations that avoid the patterns $v$ and $w$
simultaneously, begin with the pattern $p$ and end with the
pattern $q$. It is easy to see that if $k\geq 3$ or $\ell \geq 3$
then $E_{12\mn3,21\mn3}^{12\ldots k,12\ldots\ell}(x)=0$. For the
other $k$ and $\ell$, one can prove the following theorem:

\begin{thm}
We have

{\rm(i)} $E_{12\mn3,21\mn3}^{1,1}(x)=e^{x+x^2/2}-1$.

{\rm(ii)} $E_{12\mn3,21\mn3}^{1,12}(x)=e^{x+x^2/2}\left(
1-\int_0^x e^{-t-t^2/2}dt\right)-1$.

{\rm(iii)} $E_{12\mn3,21\mn3}^{12,1}(x)=\int_0^x te^{t+t^2/2} dt$.

{\rm(iv)}
$E_{12\mn3,21\mn3}^{12,12}(x)=\frac{1}{2}x^2+\int_0^x\left[e^{t+t^2/2}\left(1-\int_0^t
e^{-r-r^2/2}dr\right)-1\right]dt$.
\end{thm}

The second direction is a consideration of permutations in $S_n$
containing a pattern $v$ exactly $r$ times, beginning with some
pattern and ending with another pattern. For example, suppose that
$v=12\mn3$, $r=1$, $p=1\ldots k$, $q=1$, and $N_{v;r}^{p,q}(n)$
denotes the number of $n$-permutations that contain the pattern
$v$ exactly $r$ times, begin with the pattern $p$, and end with
the pattern $q$. It is easy to see that the only interesting case
is $1\leq k\leq 3$, since otherwise $N_{12\mn3;1}^{12\ldots
k,1}(n)=0$. Moreover, one can prove the following theorem:

\begin{thm}
Let $F_n$ denote the number of $n$-permutations containing
$12\mn3$ exactly once. Then, for all $n\geq 3$,
$$\begin{array}{l}
N_{12\mn3;1}^{1,1}(n)=F_n
N_{12\mn3;1}^{12,1}(n)=(n-1)F_{n-1}+(n-2)B_{n-2},\\
N_{12\mn3;1}^{123,1}(n)=(n-2)B_{n-3},
\end{array}$$
where $B_n$ is the $n$th Bell number, and $F_n$ is given
by~\cite[Corollarly~13]{ClaesMans2}.
\end{thm}

\end{document}